
\input amstex
\documentstyle{amsppt}

\loadbold
\loadeurm
\loadeurb
\loadeusm

\magnification=\magstep1
%
%
\hsize=6.7truein
\vsize=9.5truein
\hcorrection{-0.1truein}
\vcorrection{-0.2truein}

\font \smallrm=cmr10 at 9pt

\font \smallsl=cmsl10 at 9pt
\font \sc=cmcsc10

\baselineskip=14pt


\def \? {{\bf (???)}}
\def \loongrightarrow {\relbar\joinrel\relbar\joinrel\rightarrow}

\def \llongtwoheadrightarrow
 {\relbar\joinrel\relbar\joinrel\relbar\joinrel\twoheadrightarrow}
\def \longtwoheadrightarrow
 {\relbar\joinrel\relbar\joinrel\twoheadrightarrow}
\def \longhookrightarrow {\lhook\joinrel\relbar\joinrel\rightarrow}
\def \llonghookrightarrow
 {\lhook\joinrel\relbar\joinrel\relbar\joinrel\relbar\joinrel\rightarrow}

\def \gerg {\frak g}

\def \gerU {{\frak U}}
\def \gerF {{\frak F}}
\def \calU {{\Cal U}}
\def \calF {{\Cal F}}
\def \gersl {\frak{sl}}
\def \gergl {\frak{gl}}
\def \H {\hbox{\bf H}}
\def \gerH {{\frak H}}

\def \uqg {U_q(\gerg)}
\def \geruqg {\gerU_q(\gerg)}
\def \caluqg {\calU_q(\gerg)}

\def \uqgl {U_q(\gergl_2)}
\def \geruqgl {\gerU_q(\gergl_2)}
\def \caluqgl {\calU_q(\gergl_2)}
\def \geruegl {\gerU_{\,\varepsilon}(\gergl_2)}
\def \caluegl {\calU_\varepsilon(\gergl_2)}

\def \caluunogl {\calU_1\!(\gergl_2)}
\def \uqsl {U_q(\gersl_2)}
\def \geruqsl {\gerU_q(\gersl_2)}
\def \caluqsl {\calU_q(\gersl_2)}
\def \geruesl {\gerU_{\,\varepsilon}(\gersl_2)}

\def \caluunosl {\calU_1\!(\gersl_2)}

\def \fqm {F_q[M_2]}
\def \gerfqm {\gerF_q[M_2]}
\def \calfqm {\calF_q[M_2]}
\def \gerfunom {\gerF_1[M_2]}

\def \calfem {\calF_\varepsilon[M_2]}
\def \fqg {F_q[G]}
\def \gerfqg {\gerF_q[G]}
\def \calfqg {\calF_q[G]}

\def \fqgl {F_q[{GL}_2]}
\def \gerfqgl {\gerF_q[{GL}_2]}
\def \calfqgl {\calF_q[{GL}_2]}

\def \calfegl {\calF_\varepsilon[{GL}_2]}

\def \calfunogl {\calF_1[{GL}_2]}
\def \fqsl {F_q[{SL}_2]}
\def \gerfqsl {\gerF_q[{SL}_2]}
\def \calfqsl {\calF_q[{SL}_2]}

\def \calfesl {\calF_\varepsilon[{SL}_2]}

\def \calfunosl {\calF_1[{SL}_2]}

\def \uqgs {U_q(\gerg^*)}
\def \geruqgs {\gerU_q(\gerg^*)}
\def \caluqgs {\calU_q(\gerg^*)}

\def \uqgls {U_q\big({\gergl_2}^{\!*}\big)}
\def \geruqgls {\gerU_q\big({\gergl_2}^{\!*}\big)}
\def \caluqgls {\calU_q\big({\gergl_2}^{\!*}\big)}
\def \geruegls {\gerU_{\,\varepsilon}\big({\gergl_2}^{\!*}\big)}

\def \geruunogls {\gerU_1\!\big({\gergl_2}^{\!*}\big)}

\def \uqsls {U_q\big({\gersl_2}^{\!*}\big)}
\def \geruqsls {\gerU_q\big({\gersl_2}^{\!*}\big)}
\def \caluqsls {\calU_q\big({\gersl_2}^{\!*}\big)}
\def \geruesls {\gerU_{\,\varepsilon}\big({\gersl_2}^{\!*}\big)}

\def \geruunosls {\gerU_1\!\big({\gersl_2}^{\!*}\big)}

\def \gerFr {\frak{Fr}}
\def \calFr {\Cal{F}r}

\def \e {\text{\rm e}}
\def \f {\text{\rm f}}
\def \g {\text{\rm g}}
\def \h {\text{\rm h}}

\def \b {\hbox{\bf b}}
\def \c {\hbox{\bf c}}

\def \ebar {\overline{E}}
\def \fbar {\overline{F}}

\def \N {\Bbb N}
\def \Z {\Bbb Z}
\def \C {\Bbb C}
\def \Q {\Bbb Q}

\def \Qq {\Q (q)}
\def \Zqqm {\Z \!\left[ q, q^{-1}\right]}
\def \Qqqm {\Q \!\left[ q, q^{-1}\right]}
\def \Zeps {\Z_\varepsilon}
\def \Qeps {\Q_\varepsilon}

\document


\topmatter

{\ }

\vskip-41pt

\hfill   {\smallrm To appear in
{\smallsl Communications in Algebra\/} }
%
%
\hskip19pt   {\ }

\vskip31pt

\title
  $ \boldkey{F}_\boldkey{q} \boldkey{[}
{\boldkey{M}}_{\hskip0,3pt\boldkey{2}} \boldkey{]} \, $,
$ \boldkey{F}_\boldkey{q} \boldkey{[}
{\boldkey{G}\boldkey{L}}_{\hskip0,7pt\boldkey{2}} \boldkey{]} $
AND  $ \boldkey{F}_\boldkey{q} \boldkey{[}
{\boldkey{S}\boldkey{L}}_{\hskip0,7pt\boldkey{2}}
\boldkey{]} $  \\
  AS QUANTIZED HYPERALGEBRAS
\endtitle

\author
       Fabio Gavarini${}^\dagger$,  \ Zoran Raki\'{c}$^{\,\ddagger}$
\endauthor

\leftheadtext{ Fabio Gavarini{\,}, \ \  Zoran Raki\'{c} }
\rightheadtext{ $ F_q[M_2] \, $,  $ F_q[{GL\,}_2] $  and
$ F_q[{SL\,}_2] $  as quantized hyperalgebras }

\affil
  $ {}^\dagger $Universit\`a di Roma ``Tor Vergata'' ---
Dipartimento di Matematica  \\
  Via della Ricerca Scientifica 1, I-00133 Roma --- ITALY  \\
                             \\
  \hbox{ $ {}^\ddagger $Univerzitet u Beogradu ---
Matemati\v{c}ki fakultet }  \\
  \hbox{ Studentski trg 16, 11000 Beograd --- SERBIA }  \\
\endaffil

\address\hskip-\parindent
  Fabio Gavarini  \newline
  \indent   Universit\`a degli Studi di Roma ``Tor Vergata''
---   Dipartimento di Matematica  \newline
  \indent   Via della Ricerca Scientifica 1, I-00133 Roma, ITALY
---   gavarini\@{}mat.uniroma2.it  \newline
  {}  \newline
  Zoran Raki\'{c}  \newline
  \indent   Univerzitet u Beogradu
---   Matemati\v{c}ki fakultet  \newline
  \indent   Studentski trg 16, 11000 Beograd, SERBIA
---   zrakic\@{}matf.bg.ac.yu
\endaddress

\abstract
   Within the quantum function algebra  $ F_q[{SL}_2] $,  we study
the subset  $ \calfqsl $   --- introduced in [Ga1] ---   of all
elements of  $ \fqsl $  which are  $ \Zqqm $--valued  when paired
with  $ \caluqsl \, $,  the unrestricted  $ \Zqqm $--integral  form
of  $ \uqsl $  introduced by De Concini, Kac and Procesi.  In
particular we yield a presentation of it by generators and relations,
and a nice  $ \Zqqm $--spanning  set (of PBW type).  Moreover, we give
a direct proof that  $ \calfqsl $  is a Hopf subalgebra of  $ \fqsl $,
and that  $ \, \calfqsl\Big|_{q=1} \!\!\! \cong U_\Z({\gersl_2}^{\!*})
\, $.  We describe explicitly its specializations at roots of 1, say
$ \varepsilon $,  and the associated quantum Frobenius (epi)morphism
(also introduced in [Ga1]) from  $ \calfesl \, $  to  $ \, \calfunosl
\cong U_\Z({\gersl_2}^{\!*}) \, $.  The same analysis is done for
$ \, \calfqgl \, $,  \, with similar results, and also (as a key,
intermediate step) for  $ \, \calfqm \, $.
\endabstract

\endtopmatter

%
%

\footnote""{Keywords: \ {\sl Hopf Algebras, Quantum Groups}.}

\footnote""{ 2000 {\it Mathematics Subject Classification:} \
Primary 16W30, 17B37; Secondary 81R50. }

\footnote""{\hskip-11pt   Supported by a cooperation agreement between
the Department of Mathematics of the University of Rome ``Tor
Vergata'' and the Faculty of Mathematics of the University of
Belgrade (2000--2003)   ---
%
%
 Z.~Raki\' c
was also supported by the Serbian Ministry of Science, contract No.~144032D.}

\vskip-5pt

\centerline {\bf Introduction }

\vskip13pt

   Let  $ G $  be a semisimple, connected, simply connected affine
algebraic group over  $ \C \, $, and  $ \gerg $  its tangent Lie
algebra.  Let  $ \uqg $  be the Drinfeld-Jimbo quantum group over
$ \gerg \, $,  \, defined over the field  $ \Q(q) \, $,  \, where
$ q $  is an indeterminate.  There exist two integral forms of
$ \uqg $  over  $ \Zqqm $,  the restricted one, say  $ \geruqg
\, $,  and the unrestricted one, say  $ \caluqg $   --- see
[CP] and references therein.  Both of them bear so called
``quantum Frobenius morphisms'', namely Hopf algebra morphisms
linking their specialisations at 1 with their specialisations
at roots of 1.  In particular,  $ \geruqg $  for  $ \, q
\rightarrow 1 \, $  specializes to  $ U_\Z(\gerg) \, $,
the Kostant  $ \Z $--form  of  $ U(\gerg) \, $;  \, so
$ \gerg $  becomes a Lie bialgebra, and  $ G $  a Poisson
group.  Also,  $ \caluqg $  for  $ \, q \! \rightarrow \! 1
\, $  specializes to  $ F_\Z[G^*] \, $,  \, a  $ \Z $--form
of the function algebra on a Poisson group  $ G^* $  dual to
$ G \, $.
                                                   \par
   Dually, one constructs a Hopf algebra  $ \fqg $  of matrix
coefficients of  $ \uqg \, $.  It has two  $ \Zqqm $--forms,  say
$ \gerfqg $  and  $ \calfqg \, $,  \, defined to be the subset of
$ \fqg $  of all  $ \Zqqm $--valued  functions on  $ \geruqg \, $,
respectively on  $ \caluqg \, $.  At  $ \, q = 1 \, $,  $ \, \gerfqg $
specializes to  $ F_\Z[G] \, $,
%
%
while  $ \calfqg $  specializes to  $ U_\Z\big(\gerg^*\big) \, $,  \, a
Kostant-like  $ \Z $--form  of  $ U\big(\gerg^*\big) \, $   --- cf.~[Ga1]
for details.  Moreover, both  $ \gerfqg $  and  $ \calfqg $  bear quantum
Frobenius morphisms (relating their specialisations at 1 with those at
roots of 1), which are dual to those of  $ \geruqg $  and  $ \caluqg \, $.
                                                   \par
   The aim of this paper is to describe  $ \, \calfqg \, $,  its
specializations at roots of 1 and its quantum Frobenius morphisms
when  $ \, G = SL_2 \, $.  Moreover, as the construction of these
objects makes sense for  $ \, G = GL_2 \, $  and  $ \, G = 
M_2 := \text{\it Mat}_2 \, $  as well, we find similar results for
%
%
them.
                                                   \par
   By [Ga1],  $ \calfqm $  should resemble  $ \geruqgl \, $.  Indeed, this is the
case:  $ \calfqm $  is generated by quantum divided powers and quantum binomial
coefficients, a PBW-like theorem hold for  $ \calfqm \, $,  \, and the quantum
Frobenius morphisms are given by an  ``$ \ell $--th  root operation'',  if
$ \ell $  is the order of the root of unity.  Similar (weaker) results hold
for  $ \calfqgl $  and  $ \calfqsl \, $.
%
%
                                                   \par
   The general case of  $ \, M_n := \text{\it Mat}_n \, $,  $ {GL}_n $
and  $ {SL}_n $  is studied in [GR2], exploiting the same key ideas
already developed here and the present results for  $ \, n = 2 \, $.
 \vskip2pt
   {\it Warning:} \, an expanded, more detailed version of this paper is
available on line, cf.~[GR1]; the quotations in [GR2] about the present
work refer in fact to [GR1].

\vskip13pt

   \centerline{\sc dedicatory}
 \vskip-1pt
  {\sl This work grew out of a cooperation supported by an
official agreement between the Department of Mathematics of
the University of Rome ``Tor Vergata'' and the Faculty of
Mathematics of the University of Belgrade in the period
2000--2003.  Such agreement was the outcome of a common
wish of peaceful, fruitful partnership, as an answer to
the military aggression of NATO countries to the Federal
Republic of Yugoslavia, which started in springs of 1999.
This paper is dedicated to the memory of all victims of
that war. }

\vskip1,1truecm

\centerline {\bf \S\; 1 \ Geometrical background and  $ q $--numbers }

\vskip13pt

   {\bf 1.1 Poisson structures on linear groups.}  Let  $ \, \gerg :=
\gergl_2({\Bbb Q}) \, $,  \, with its basis given by the elementary
matrices  $ \, e := \Big({{0 \ 1} \atop {0 \ 0}}\Big) \, $,  $ \,
g_1 := \Big({{1 \ 0} \atop {0 \ 0}}\Big) \, $,  $ \, g_2
:= \Big({{0 \ 0} \atop {0 \ 1}}\Big) \, $,  $ \, f :=
\Big({{0 \ 0} \atop {1 \ 0}}\Big) \, $.  Then  $ \gerg $
has a natural structure of Lie algebra, and a Lie cobracket
is defined on it by  $ \, \delta(e) = h \otimes e - e \otimes
h \, $,  $ \, \delta(g_k) = 0 \, $  (for $ \, k = 1, 2 \, $),  $ \,
\delta(f) = h \otimes f - f \otimes h \, $,  \, where  $ \, h :=
g_1 - g_2 \, $;  \, this makes  $ \gerg $  into a  {\sl Lie
bialgebra}.  It follows that  $ \, U(\gerg) \, $  is naturally
a co-Poisson Hopf algebra, whose co-Poisson bracket is the extension
of the Lie cobracket of  $ \gerg \, $.
%
%
 Finally, Kostant's  $ \Z $--integral  form of  $ U(\gerg) $   ---
called also  {\sl hyperalgebra\/}  in literature ---   is the unital
$ \Z $--subalgebra  $ U_\Z(\gerg) $  of  $ U(\gerg) $  generated
by the ``divided powers''  $ \, f^{(n)} \, $,  $ \, e^{(n)} \, $  and
the binomial coefficients  $ \, \Big(\! {g_k \atop n} \!\Big) \, $
(for  $ \, k = 1, 2 \, $,  and  $ \, n \in \N \, $),  where
%
%
we use notation  $ \, x^{(n)} := x^n \big/ n! \, $  and  $ \, \Big(\!
{t \atop n} \!\Big) := {{\, t (t-1) \cdots (t-n+1) \,} \over {\, n! \,}}
\, $.  Again, this is a co-Poisson Hopf  $ \Z $--algebra;  it is free
as a  $ \Z $--module,  with PBW-like  $ \Z $--basis  the set of ordered
monomials  $ \; \Big\{ e^{(\eta)} \Big(\! {g_1 \atop \gamma_1} \!\Big)
\Big(\! {g_2 \atop \gamma_2} \!\Big) f^{(\varphi)} \,\Big|\, \eta,
\gamma_1, \gamma_2, \varphi \in \N \,\Big\} \, $; see e.g.~[Hu],
Ch.~VII.
                                           \par
   A similar description holds for  $ \, \gerg := \gersl_2(\Q) \, $,
\, taking  $ h $  instead of  $ g_1 $  and  $ g_2 \, $.  The Kostant's
$ \Z $--form  $ \, U_\Z\big(\gersl_2\big) \, $  of  $ U \big( \gersl_2
(\Q) \big) $  is generated
%
%
 as above but for replacing the  $ g_k $'s  with  $ h \, $.
Then  $ \, U_\Z\big(\gersl_2\big) \, $  is a co-Poisson Hopf subalgebra
of  $ \, U_\Z\big(\gergl_2\big) \, $,  \, free as a  $ \Z $--module
with PBW  $ \Z $--basis  as above but with  $ h $  instead of the
$ g_k $'s.  Finally,  $ \, \gersl_2(\Q) \, $  is a Lie sub-bialgebra
of  $ \gergl_2(\Q) \, $,  \, and the embedding  $ \; \gersl_2
\lhook\joinrel\relbar\joinrel\rightarrow \gergl_2
\; $  is a section of the natural Lie bialgebra epimorphism  $ \; \gergl_2
\relbar\joinrel\twoheadrightarrow \gersl_2 \; $.
%
%
                                            \par
   As  $ \gergl_2(\Q) $  is a Lie bialgebra, by general theory  $ \,
G := {GL}_2(\Q) \, $  is then a  {\sl Poisson\/}  group.
%
%
 Explicitly,  the
algebra  $ F[G] $  of regular functions on  $ G $  is the unital
associative commutative  $ \Q $--algebra  with generators  $ \,
\bar{a} $,  $ \bar{b} $,  $ \bar{c} $,  $ \bar{d} $  and  $ D^{-1} $,
where  $ \, D := \text{\sl det}\,\Big({{\bar{a} \ \bar{b}} \atop
{\bar{c} \ \bar{d}}}\Big) \, $  is the determinant.  The group
structure on  $ G $  yields on  $ F[G] $  the natural Hopf structure
given by matrix product,
%
%
 while the
Poisson structure is given by
 \vskip-17pt
  $$  \big\{\bar{a}, \bar{b}\big\} = \bar{b} \, \bar{a} \, ,  \hskip8pt
\{\bar{a}, \bar{c}\} = \bar{c} \, \bar{a} \, ,  \hskip8pt  \big\{\bar{b},
\bar{c}\big\} = 0 \, ,  \hskip8pt  \big\{\bar{d}, \bar{b}\big\} = -
\bar{b} \, \bar{d} \, ,  \hskip8pt  \big\{\bar{d}, \bar{c}\big\} = -
\bar{c} \, \bar{d} \, ,  \hskip8pt  \big\{\bar{a}, \bar{d}\,\big\} =
2 \, \bar{b} \, \bar{c} \, .  $$
 \vskip-7pt
   \indent   We shall consider also the Poisson group-scheme  $ G_\Z $
associated to  $ {GL}_2 $,  \, for which a like analysis applies: in
particular, its function algebra  $ F[G_\Z] $  is a Poisson Hopf
$ \Z $--algebra  with the same presentation as  $ F[G] $  but over
the ring  $ \Z \, $.
 \vskip2pt
   Similar constructions hold for  $ {SL}_2(\Q) $  and the associated group-scheme
(just set  $ D = 1 \, $).
 \vskip2pt
   Finally, the subalgebra of  $ F\big[({GL}_2)_{\Z}\big] $  generated
by the  $ \bar{a} $,  $ \bar{b} $,  $ \bar{c} $  and  $ \bar{d} $
is a Poisson subbialgebra of  $ F\big[({GL}_2)_{\Z}\big] \, $:  \,
indeed, it is the algebra  $ F\big[(M_2)_{\Z}\big] $  of regular
functions of the  $ \Z $-scheme  associated to the Poisson algebraic
monoid  $ \, M_2 \, $  of all  $ (2 \times 2) $--matrices.

\vskip11pt

   {\bf 1.2 Dual Lie bialgebras and dual Poisson groups.}  By general
theory, if  $ \, \gerg := \gergl_2(\Q) \, $  bears a Lie bialgebra
structure then the dual space  $ \, \gerg^* \, $  is a Lie bialgebra on
its own.  Let  $ \, \big\{ e^*, g_1^{\,*}, g_2^{\,*}, f^* \big\} \, $
be the dual basis to the basis of elementary matrices for  $ \gerg \, $,
and let  $ \, \e := f^* \!\big/ 2 \, $,  $ \, \g_1 := g_1^{\,*} \, $,
$ \, \g_2 := g_2^{\,*} \, $,  $ \, \f := e^* \!\big/ 2 \, $;  \, then
$ \, \big\{ \e \, , \g_1 \, , \g_2 \, , \f \,\big\} $  is a basis of
$ \gerg^* $.  The Lie bracket of  $ \gerg^* $  is given by
$ \; [\, \g_1 \, , \, \g_2 \,] = 0 \, $,  $ \, [\, \g_1 \, ,
\, \f \,] = +\,\f \, $,  $ \, [\, \g_2 \, , \, \f \,] =
-\,\f \, $,  $ \, [\, \g_1 \, , \, \e \,] \, = \, +\,\e \, $,
$ \, [\, \g_2 \, , \, \e \,] = -\,\e \, $,  $ \, [\, \f \, ,
\, \e \,] = 0 \, $,  \; and its Lie cobracket by
$ \; \delta(\f\,) \! = \! (\g_1 - \g_2) \wedge \f \, $,
$ \, \delta(\g_1) \! = 4 \cdot \f \wedge \e \, $,
$ \, \delta(\g_2) \! = 4 \cdot \e \wedge \f \, $,
$ \, \delta(\e) \! = \e \wedge (\g_1 - \g_2) \, $,
\; where  $ \, x \wedge y := x \otimes y - y \otimes x \, $.
These formul{\ae}  also provide a presentation of  $ U(\gerg^*) $
as a co-Poisson Hopf algebra.  Finally, we can define the Kostant's
$ \Z $--integral  form, or hyperalgebra,  $ U_\Z(\gerg^*) $  of
$ U(\gerg^*) $  as the unital  $ \Z $--subalgebra  generated by
the divided powers  $ \, \f^{\,(n)} \, $,  $ \, \e^{(n)} \, $
and binomial coefficients  $ \Big(\! {\g_k \atop n} \!\Big) $  (for
all  $ \, n \in \N \, $  and all  $ \, k = 1, 2 \, $).  This again is
a co-Poisson Hopf  $ \Z $--algebra,  free as a  $ \Z $--module,  with
PBW-like  $ \Z $--basis  the set of ordered monomials  $ \; \Big\{
\e^{(\eta)} \Big(\! {\g_1 \atop n_1} \!\Big) \Big(\! {\g_2 \atop n_2}
\!\Big) \f^{\,(\varphi)} \;\Big|\; \eta \, , n_1 \, , n_2 \, , \varphi
\in \N \,\Big\} \, $.
                                           \par
   A like description holds for  $ \, {\gersl_2(\Q)}^* \, $:  \, indeed,
one has  $ \, {\gersl_2 (\Q)}^* \! = {\gergl_2(\Q)}^* \!\Big/ \big( \g_1
\! + \g_2 \big) \, $,  \, dually to  $ \, \gersl_2(\Q) \longhookrightarrow
\gersl_2(\Q) \, $,  hence one simply has to set  $ \, \h := \g_1 \equiv -
\g_2 \, $  in the presentation above.  All formul{\ae}  involving  $ \h $
follow from  $ \, \h \cong +\g_1 \cong -\g_2 \, \mod \big( \g_1 \! + \g_2 \big)
\, $.  In particular  $ U_\Z\big(\gersl_2^{\,*}\big) $  is the  $ \Z $--subalgebra  of  $ U\big(\gersl_2(\Q)^*\big) $  generated by divided powers and binomial
coefficients as above but taking  $ \h $  instead of the  $ \g_k $'s.  Then
$ \, U_\Z\big(\gersl_2^{\,*}\big) \, $  is a co-Poisson Hopf  $ \Z $--subalgebra
of  $ \, U_\Z\big(\gergl_2^{\,*} \big) \, $,  \, with PBW  $ \Z $--basis  as
above but with  $ \h $  instead of the  $ \g_k $'s.
                                            \par
   If  $ \, \gerg = \gergl_2 \, $,  \, a simply connected algebraic
Poisson group with tangent Lie bialgebra  $ \gerg^* $  is the subgroup
$ {}_s{G}^* $  of  $ \, G \times G \, $  made of all pairs  $ \, \big(
L, U \big) \in G \times G \, $  such that  $ L $  is lower triangular,
$ U $  is upper triangular, and their diagonals are inverse to each other.
This is a Poisson subgroup of  $ G \times G \, $;  \, its centre is  $ \,
Z := \big\{ \big(z I , z^{-1} I \big) \,\big|\, z \in \Q \setminus \{0\}
\big\} \, $,  \, hence the associated adjoint group is  $ \, {}_a{G}^* :=
{}_s{G}^*\! \big/ Z \, $.  The same construction defines Poisson group-schemes
$ {}_s{G}^*_\Z \, $  and  $ {}_a{G}^*_\Z \; $.  If  $ \, \gerg = \gersl_2 \, $
the construction of dual Poisson group-schemes  $ {}_s{G}^*_\Z \, $  and
$ {}_a{G}^*_\Z \; $  is entirely similar, just taking  $ \, G := {SL}_2 \, $
instead of  $ {GL}_2 $  in the previous recipe.

\vskip11pt

   {\bf 1.3  $ q $--numbers,  $ q $--divided  powers and  $ q $--binomial
coefficients.}  Let  $ q $  be an indeterminate.  For all  $ \, s, n \in
\N \, $,  let  $ \; {(n)}_q := {{q^n - 1} \over {q - 1}} \; (\in \Z[q])
\, $,  $ \; {(n)}_q! := \prod_{r=1}^n {(r)}_q \, $,  $ \; {\Big( {n \atop s}
\Big)}_{\!q} := {{{(n)}_q!} \over {{(s)}_q! {(n-s)}_q!} } \;\; (\in \Z[q])
\, $,  \; and  $ \; {[n]}_q := {{q^n - q^{-n}} \over {q - q^{-1}}} \; (\in
\Zqqm) \, $,  $ \; {[n]}_q! := \prod_{r=1}^n {[r]}_q \, $,  $ \; {\Big[
{n \atop s} \Big]}_q \! := {{{[n]}_q!} \over {{[s]}_q! {[n-s]}_q!} } \;\;
(\in \Zqqm) \, $.  Furthermore, we set  $ {\Big( {-n \atop s} \Big)}_{\! q}
\! := {(-1)}^s q^{- n s - {s \choose 2}} {{{(n \! - \! 1 \! + \! s)}_q!}
\over {(s)}_q! } \; (\in \Z[q]) \, $  for all  $ \, n, s \in \N \, $,
and  $ \, {(2k)}_q!! := \prod_{r=1}^k {(2r)}_q \, $,  $ \, {(2k-1)}_q!!
:= \prod_{r=1}^k {(2r-1)}_q \, $,  \, for all  $ \, k \in \N_+ \, $.
                                       \par
   If  $ A $  is any  $ \Qq $--algebra,  $ q $--divided  powers and
$ q $--binomial  coefficients are:  $ \; X^{(n)} := X^n \Big/ {[n]}_q!
\; $,  $ \; \Big(\! {{X \, ; \, c} \atop n} \!\Big) := \prod_{s=1}^n {{\,
q^{c+1-s} X - 1 \,} \over {\, q^s - 1 \,}} \; $;  \; also,
$ \; \left\{\! {{X \, ; \, c} \atop {n \, , \, r}} \right\} \, := \,
\sum_{s=0}^r q^{{s+1} \choose 2} {\Big(\! {r \atop s} \!\Big)}_{\!q}
\cdot \Big(\! {{X \, ; \, c+s} \atop {n-r}} \!\Big) \; $
\; (for every  $ \, X \in A \, $,  $ \, n, r \in \N \, $,
$ \, c \in \Z \, $).  Furthermore, if
$ \, Z \in A \, $  is  {\sl invertible\/}  we define also  $ \; \left[\!
{{Z \, ; \, c} \atop n} \right] := \prod_{s=1}^n {{\, q^{c+1-s} Z^{+1}
- q^{s-1-c} Z^{-1} \,} \over {\, q^{+s} - q^{-s} \,}} \; $  for every
 $ \, n \in \N \, $  and  $ \, c \in \Z \, $.
                                             \par
   For later use, we remark that the  $ q $--binomial  coefficients
in  $ \, X \in A \, $  satisfy the relations
 \vskip-13pt
  $$  \hbox{ $ \eqalign{
   {}  &  {\textstyle \prod_{s=1}^n} \big( q^s \! - \! 1 \big)
\textstyle{ \Big(\! {{X \, ; \; c} \atop n} \!\Big) } \, = \, {\textstyle \prod_{s=1}^n} \big( q^{1-s+c} X - 1 \big)  \cr
   \textstyle{ \Big(\! {{X \, ; \; c} \atop t} \!\Big)
\Big(\! {{X \, ; \; c-t} \atop s} \!\Big) } =  &
\; \textstyle{ \Big(\! {{t+s} \atop t} \!\Big)_{\! q}
\Big(\! {{X \, ; \; c} \atop {t+s}} \!\Big) } \; ,  \qquad
\textstyle{ \Big(\! {{X \, ; \; c+1} \atop t} \!\Big) }
- q^t \textstyle{ \Big(\! {{X \, ; \; c} \atop t} \!\Big) } =
\textstyle{ \Big(\! {{X \, ; \; c} \atop {t-1}} \!\Big) }  \cr
   \textstyle{ \Big(\! {{X \, ; \; c} \atop m} \!\Big) \Big(\! {{X \, ; \; s}
\atop n} \!\Big) } =  &  \; \textstyle{ \Big(\! {{X \, ; \; s} \atop n} \!\Big) \Big(\! {{X \, ; \; c} \atop m} \!\Big) } \; ,  \quad  \textstyle{ \Big(\!
{{X \, ; \; c} \atop t} \!\Big) } = {\textstyle \sum_{p \geq 0}^{p \leq c,t}} \, q^{(c-p)(t-p)} \textstyle{ \Big(\! {c \atop p} \!\Big)_{\! q} } \textstyle{ \Big(\! {{X \, ; \; 0} \atop {t-p}} \!\Big) }  \cr
   \textstyle{ \Big(\! {{X \, ; \; c} \atop 0} \!\Big) } = 1 \; ,  \quad  &  \textstyle{ \Big(\! {{X \, ; \; -c} \atop t} \!\Big) } =
{\textstyle \sum_{p=0}^t} \, {(-1)}^p q^{- t(c+p) + p(p+1)/2}
\textstyle{ \Big(\! {{p+c-1} \atop p} \!\Big)_{\! q} \Big(\! {{X \, ; \; 0}
\atop {t-p}} \!\Big) }  \cr
   \textstyle{ \Big(\! {{X \, ; \; c+1} \atop t} \!\Big) } -  &
\textstyle{ \Big(\! {{X \, ; \; c} \atop t} \!\Big) } = \, q^{c-t+1}
\left( 1 + (q-1) \textstyle{ \Big(\! {{X \, ; \; 0} \atop 1} \!\Big) } \right)
\textstyle{ \Big(\! {{X \, ; \; c} \atop {t-1}} \!\Big) }  \cr } $ }
\eqno (1.1)  $$
 \vskip-3pt
\noindent
 Similarly, the  $ q $--divided  powers in  $ \, X \in A \, $
satisfy the relations
 \vskip-8pt
  $$  X^{(r)} X^{(s)} \, = \textstyle{ \Big[\! {{r+s} \atop s} \!\Big]_q } X^{(r+s)} \; ,  \qquad  X^{(0)} = 1   \eqno (1.2)  $$
 \vskip-2pt
%
%
 \vskip-3pt
   \indent   Finally, let  $ \, \ell \in \N_+ \, $  be odd, set  $ \,
\Zeps := Z[q] \Big/ \big(\phi_\ell(q)\big) \, $  where  $ \, \phi_\ell(q)
\, $  is the  $ \ell $-th  cyclotomic polynomial in  $ q \, $,  \, and
let  $ \, \varepsilon := \overline{q} \, $,  \, a (formal) primitive
$ \ell $-th  root of 1 in  $ \Zeps \, $.  Similarly let  $ \, \Qeps :=
\Q[q] \Big/ \big(\phi_\ell(q)\big) \, $,  \, the field of quotients of
$ \Zeps \, $.  If  $ M $  is a module over  $ \Zqqm $  or  $ \Qqqm $
we shall set  $ \; M_\varepsilon := M \Big/ \big(\phi_\ell(q)\big) \,
M \, $,  \, which is a module over  $ \Zeps $  or over  $ \Qeps \, $.

\vskip1truecm

\centerline {\bf \S\; 2 \ Quantum groups }

\vskip13pt

  {\bf 2.1 Quantum enveloping algebras  $ \uqgl $  and  $ \uqsl \, $,  their
integral forms and specializations.} \, Let  $ \uqgl $  be the unital
$ \Qq $--algebra  with generators  $ \, F $,  $ G_1 \, $,  $ G_2 \, $,
$ {G_1}^{\! -1} $,  $ {G_2}^{\! -1} $,  $ E \, $  and relations
 \vskip-17pt
  $$  \displaylines{
   {G_i}^{\! \pm 1} {G_i}^{\! \mp 1} = 1 = {G_i}^{\! \mp 1} {G_i}^{\! \pm 1}
\,\; \big( i=1, 2 \big) \, ,  \hskip8pt  G_1 G_2 = G_2 \, G_1 \; ,  \hskip8pt
E F - F E \, = {{\; G_1 {G_2}^{\!\! -1} \! - {G_1}^{\!\! -1} G_2 \;}
\over {\; q - q^{-1} \;}}  \cr
   {G_1}^{\! \pm 1} F = q^{\mp 1} F {G_1}^{\! \pm 1} ,  \;\;
{G_2}^{\! \pm 1} F = q^{\pm 1} F {G_2}^{\! \pm 1} ,  \;\;
{G_1}^{\! \pm 1} E = q^{\pm 1} E {G_1}^{\! \pm 1} ,  \;\;
{G_2}^{\! \pm 1} E = q^{\mp 1} E {G_2}^{\! \pm 1}  \cr }  $$
%
%
 \vskip-0pt
   \indent   Moreover,  $ \uqgl $  is also a Hopf algebra with
$ \; \Delta(F) = F \otimes {G_1}^{\! -1} G_2 + 1 \otimes F \, $,
$ \epsilon(F) = 0 \; $,  $ \; S(F) = - F G_1 {G_2}^{\! -1} \; $,
$ \; \Delta\big({G_i}^{\!\pm 1}\big) = {G_i}^{\!\pm 1} \otimes
{G_i}^{\!\pm 1} \; $,  $ \; \epsilon\big({G_i}^{\!\pm 1}\big) =
1 \; $,  $ \; S\big({G_i}^{\! \pm 1}\big) = {G_i}^{\!\mp 1} \; $
(for  $ \, i = 1, 2 \, $),  $ \; \Delta(E) = E \otimes 1 + G_1
{G_2}^{\! -1} \otimes E \; $,  $ \; \epsilon(E) = 0 \; $,  $ \;
S(E) = - {G_1}^{\! -1} G_2 \, E \; $.
                                             \par
   Now let  $ \; K^{\pm 1} := G_1^{\pm 1} G_2^{\mp 1} \, $.  The
unital  $ \Qq $--subalgebra  of  $ \uqgl $  generated by  $ \, F $,
$ K^{\pm 1} $  and  $ E $  is the well-known Drinfeld-Jimbo's quantum
algebra  $ \uqsl \, $.  From the presentation of  $ \uqgl $  one argues
one of  $ \uqsl $  too, and also sees that the latter is a Hopf
subalgebra of the former; indeed, it is also a quotient via  $ \,
F \mapsto F \, $,  $ \, G_1 \mapsto K \, $,  $ \, G_2 \mapsto K^{-1}
\, $,  $ \, E \mapsto E \, $.
                                           \par
   The quantum version of the PBW theorem for  $ \uqgl $  claims that
the set of ordered monomials  $ \; B^g := \left\{\, E^\eta \, {G_1}^{\!
\gamma_1} \, {G_2}^{\! \gamma_2} \, F^\varphi \,\big|\, \eta \, ,
\varphi \in \N \, , \, \gamma_1, \gamma_2 \in \Z \,\right\} \; $  is
a  $ \Qq $--basis  of  $ \uqgl \, $.  Similarly, the set  $ \; B^s :=
\left\{\, E^\eta \, K^\kappa \, F^\varphi \,\big|\, \eta \, , \varphi
\in \N \, , \, \kappa \in \Z \,\right\} \; $  is a  $ \Qq $--basis
of  $ \uqsl \, $.
 \vskip1pt
  As to integral forms, let  $ \geruqgl $  be the unital  $ \Zqqm $--subalgebra
of  $ \uqgl $  generated by  $ \; F^{(m)} $,  $ G_1^{\pm 1} $,  $ \Big(\!
{{G_1 \, ; \, c} \atop m} \!\Big) $,  $ G_2^{\pm 1} $,  $ \Big(\! {{G_2 \, ; \, c}
\atop m} \!\Big) $,  $ E^{(m)} $,  \, for all  $ \, m \in \N \, $  and  $ \, c
\in \Z \, $.  This  $ \geruqgl $  is a  $ \Zqqm $--integral  form of  $ \uqgl $
as a Hopf algebra, and specializes to  $ U_\Z(\gergl_2) $  for  $ \, q \mapsto 1
\, $,  \, that is  $ \; \geruqgl \Big/ (q-1) \, \geruqgl \cong U_\Z(\gergl_2) \; $  as co-Poisson Hopf algebras; therefore we call  $ \, \geruqgl \, $  a  {\it quantum (\/{\rm or}  quantized) hyperalgebra}.
%
%
 Finally, for every root of 1, say  $ \varepsilon $,  of odd order  $ \ell \, $,
a  {\it quantum Frobenius morphism\/}  $ \; {\frak{Fr}}_{\gergl_n}^{\,\Z} \,
\colon \, \geruegl \relbar\joinrel\relbar\joinrel\twoheadrightarrow \Zeps
\otimes_\Z U_\Z(\gergl_n) \; $  exists (a Hopf algebra epimorphism): the
left-hand side is  $ \, \geruegl := {(\geruqgl)}_\varepsilon \, $,  \,
and  $ {\frak{Fr}}_{\gergl_n}^{\,\Z} $  is defined on generators by
``dividing out by  $ \ell \, $''  the order of each quantum divided
power and each quantum binomial coefficient, if this makes sense,
and mapping to zero otherwise.
                                           \par
   Similarly, one defines the integral form of  $ \uqsl \, $,  say
$ \geruqsl \, $,  \, replacing the  $ G_k^{\pm 1} $'s  by  $ K^{\pm 1} $,
and the  $ \Big(\! {{G_k \, ; \, c} \atop m} \!\Big) $'s  by the  $ \Big(\!
{{K \, ; \, c} \atop m} \!\Big) $'s;  \, totally similar results then hold
%
%
(see [DL], [Ga1]).
                                           \par
   As to  {\sl unrestricted\/}  integral forms and their specializations, first
set  $ \, \overline{X} := \big(q - q^{-1}\big) X \, $  as notation.
%
 We
define  $ \caluqgl $  to be the unital
$ \Zqqm $--subalgebra  of  $ \uqgl $  generated by  $ \; \big\{\,
\fbar \, , \, G_1^{\pm 1} , \, G_2^{\pm 1} , \, \ebar \,\big\} \; $.
   From the presentation of  $ \uqgl $  one argues a presentation for
$ \caluqgl $  as well, and then sees that the latter is a Hopf subalgebra
of the former.  Moreover,  $ \caluqgl $  is a free  $ \Zqqm $--module  with
basis  $ \; \Cal{B}^g := \Big\{\, \fbar^{\,\varphi} G_1^{\gamma_1} G_2^{\gamma_2}
\ebar^{\,\eta} \,\Big|\, \varphi, \gamma_1, \gamma_2, \eta \in \N \Big\} \, $.
Note that  $ \caluqgl $  is another  $ \Zqqm $--integral  form of  $ \uqgl \, $,
as a Hopf algebra, in that it is a Hopf  $ \Zqqm $--subalgebra  such that  $ \,
\Qq \otimes_{\Z[q,q^{-1}]} \caluqgl \cong \uqgl \, $.
                                        \par
   Adapting results in [DP], [Ga1] and [Ga3--4], one has that  $ \caluqgl $  is
a  {\sl quantization\/}  of  $ \, F \big[{({{}_s{GL}_2}^{\!*})}_\Z\big] \, $,
\,  {i.e.}  $ \; \caluunogl := {(\caluqgl)}_1 \cong F \big[ {({}_s
{{GL}_2}^{\!*})}_\Z \big] \; $  as Poisson Hopf algebras, where on
left-hand side we consider the standard Poisson structure inherited
from  $ \caluqgl \, $.  Finally, let  $ \ell $  and  $ \varepsilon $
be as in \S\; 1.3.  Set  $ \; \caluegl := {(\caluqgl)}_\varepsilon
\, $:  \, then there is a Hopf algebra monomorphism  $ \;
{\calF{}r}_{\gergl_2}^{\,\Z} \, \colon \, F \big[ {({{}_s
{GL}_2}^{\!*})}_\Z \big] \cong \, \Zeps \otimes_\Z \caluunogl
\, \llonghookrightarrow \, \caluegl \; $  given by  $ \; \fbar
\,\big\vert_{q=1} \! \mapsto \fbar^{\,\ell}\Big\vert_{q=\varepsilon}
\; $,  $ \;  {G_k}^{\! \pm 1} \big\vert_{q=1} \!\! \mapsto {G_k}^{\!
\pm \ell}\big\vert_{q=\varepsilon} \; $,  $ \; \ebar\,\big\vert_{q=1}
\! \mapsto \ebar^{\,\ell}\Big\vert_{q=\varepsilon} \; $  ($ \, k = 1,
2 \, $).  This is the  {\it quantum Frobenius morphism\/}  for  $ \,
{{}_s{GL}_2}^{\!*} \; $.
%
%
                                        \par
   Again, the same constructions can be done with  $ \gersl_2 $  too.
One defines the unrestricted  $ \Zqqm $--integral  form  $ \caluqsl $
of  $ \uqsl \, $,  simply following the recipe above but replacing the
$ G_k^{\pm 1} $'s  with  $ K^{\pm 1} \, $.  Then similar results to
those for  $ \caluqgl $  hold for  $ \uqsl $  as well, {e.g.}
%
%
 $ \caluqsl $  is a  {\sl quantization\/}  of  $ \, F \big[
{({{}_a{SL}_2}^{\!*})}_\Z \big] \, $,  \, that is  $ \; \caluunosl
:= {(\caluqsl)}_1 \cong F \big[{( {{}_a {SL}_2}^{\!*})}_\Z \big] \; $
as Poisson Hopf algebras (like above).
%
%
 See [Ga1] and references therein for further details.
                                        \par
   The embedding  $ \; \uqsl \longhookrightarrow \uqgl \; $  restricts
to Hopf embeddings  $ \; \geruqsl \longhookrightarrow \geruqgl \, $  and
$ \; \caluqsl \longhookrightarrow \caluqgl \, $.  The specializations of
the latter ones at  $ \, q = \varepsilon \, $  and at  $ \, q = 1 \, $
are compatible (in the obvious sense) with the quantum Frobenius
morphisms.

\vskip11pt

  {\bf 2.2 Quantum function algebras  $ \fqm $,  $ \fqgl $  and  $ \fqsl \, $,
their integral forms and specializations.} \, Let  $ \fqm $  be the well-known
quantum function algebra over  $ M_2 $  introduced by Manin.  Namely,  $ \fqm $
is the unital associative  $ \Qq $--algebra  with generators  $\; a \, $,
$ \, b \, $,  $ \, c \, $,  $ \, d \; $  and relations
 \vskip-17pt
  $$  a b = q \, b a \; ,  \hskip7pt  a c = q \, c a \; ,
\hskip7pt  b d = q \, d b \; ,  \hskip7pt  c d = q \, d c \; ,
\hskip7pt  b c = c b \; ,  \hskip7pt  a d - d a = \big( q - q^{-1} \big) \, b c \; .  $$
 \vskip-5pt
\noindent
 This is also a  $ \Qq $--bialgebra  (yet not a Hopf algebra), with
coalgebra structure given by
 \vskip-17pt
  $$  \displaylines{
   \Delta(a) = a \otimes a + b \otimes c \; ,  \quad  \epsilon(a) = 1
\; ,  \qquad  \Delta(b) = a \otimes b + b \otimes d \; ,  \quad
\epsilon(b) = 0 \;\, \phantom{.}  \cr
  \Delta(c) = c \otimes a + d \otimes c \; ,  \quad  \epsilon(c) = 0
\; ,  \qquad  \Delta(d) = c \otimes b + d \otimes d \; ,  \quad
\epsilon(d) = 1 \; .  \cr }  $$
 \vskip-5pt
   \indent  In particular, the quantum determinant  $ \; D_q := a d - q \, b c
\in \fqm \; $  is central and group-like in the bialgebra  $ \fqm $.
   Finally, it follows from definitions that  $ \fqm $  admits as  $ \Qq $--basis
the set of ordered monomials  $ \; B_{M_2} := \Big\{\, b^\beta a^\alpha d^\delta
c^\kappa \,\Big|\, \beta, \alpha, \delta, \kappa \in \N \,\Big\} \, $.

\vskip2pt

   We define  $ \gerfqm $  to be the unital  $ \Zqqm $--subalgebra  of
$ \fqm $  generated by  $ \, a $,  $ b $,  $ c \, $  and  $ d \, $;  \,
this in fact is a  {\sl sub-bialgebra},  and admits the same presentation
as  $ \fqm $  but over  $ \Zqqm \, $.  It follows that  $ B_{M_2} $  is also
a  $ \Zqqm $--basis  of  $ \gerfqm \, $,  hence  $ \gerfqm $  is a  {\sl
$ \Zqqm $--integral  form of  $ \fqm \, $}.  From its presentation one sees
that  $ \gerfqm $  is a  {\sl quantization\/}  of  $ F\big[(M_2)_{\Z}\big] $,  i.e.~$ \; \gerfunom := {(\gerfqm)}_1 \cong F\big[(M_2)_\Z\big] \; $  as Poisson bialgebras (where on left-hand side we consider the standard Poisson structure
inherited from  $ \gerfqm \, $);  using notation of \S\; 1.1, the
isomorphism is given by  $ \, x{\big|}_{q=1} \cong \bar{x} \, $
for all  $ \, x \in \{a,b, c,d\,\} \, $.

 \vskip2pt

   Let  $ \; \fqgl := \big(\fqm\big)\big[{D_q}^{\!-1}\big] \, $,  \; the
extension of  $ \fqm $  by a formal inverse to  $ D_q \, $.  Then  $ \fqgl $
is the unital associative  $ \Qq $--algebra  with generators  $ \, a $,  $ b $,
$ c $,  $ d \, $  and  $ {D_q}^{\! -1} $,  \, and relations like for  $ \fqm $
plus those saying that  $ {D_q}^{\! -1} $  is central and inverse to  $ D_q \, $.
This  $ \fqgl $  is a Hopf algebra, whose coproduct and counit on the generators is
given by the formul{\ae}  in \S\; 2.2 and by those saying that  $ {D_q}^{\! -1} $
is group-like plus  $ \, S(a) = d \, {D_q}^{\! -1} \, $,  $ \, S(b) = - q^{-1} b
\, {D_q}^{\! -1} \, $,  $ \, S(c) = - q^{+1} c \, {D_q}^{\! -1} \, $,  $ \, S(d)
= a \, {D_q}^{\! -1} \, $  and  $ \, S\big({D_q}^{\! -1}\big) = a b -q \, b c
\equiv D_q \; $.
                                               \par
   It follows that  $ \fqgl $  is  $ \Qq $--spanned  by  $ \, B_{GL_2}
:= \Big\{\, b^\beta a^\alpha d^\delta c^\kappa D_q^{-n} \,\Big|\, \beta,
\alpha, \delta, \kappa, n \in \N\,\} \; $.

 \vskip2pt

   Let  $ \gerfqgl $  be the unital  $ \Zqqm $--subalgebra  of  $ \fqgl $
generated by  $ \, a $,  $ b $,  $ c $,  $ d \, $  and  $ \, {D_q}^{\! -1}
\, $  (note that  $ \, D_q \in \gerfqm \, $).  This in fact is a Hopf
$ \Zqqm $--subalgebra, and admits the same presentation as  $ \fqgl $  but
over  $ \Zqqm $.  Then  $ B_{GL_2} $  is also a  $ \Zqqm $--spanning  set of
$ \gerfqgl \, $,  \, hence  $ \gerfqgl $  is a  {\sl  $ \Zqqm $--integral
form of  $ \uqgl \, $}.  Also,  $ \gerfqgl $  is a  {\sl quantization\/}
of  $ F\big[({GL}_2)_\Z\big] $
%
%
 as a Poisson Hopf algebra, with  $ \, D_q^{\pm 1} {\big|}_{q=1} \cong
D^{\pm 1} \; $.

 \vskip2pt

   Let  $ \fqsl $  be the quotient  $ \; \fqsl := \fqgl \Big/ \! \big( D_q
- 1\big) \, \cong \, \fqm \Big/ \! \big( D_q - 1 \big) \; $  where  $ \,
\big(D_q - 1\big) \, $  is the two-sided ideal of $ \fqgl $  or of
$ \fqm $  generated by the central element  $ \, D_q - 1 \; $.  This
is a Hopf ideal of  $ \fqgl $,  so  $ \fqsl $  is a Hopf algebra too:
it admits the like presentation as  $ \fqm $  or  $ \fqgl $  but with
the additional relation  $ \; D_q - 1 = 0\; $.  Moreover,  $ \fqsl $  has
the Hopf structure given as for  $ \fqgl $  but setting  $ \; {D_q}^{\! -1}
= 1 \; $.  It also follows that  $ \fqsl $  admits  $ \; B_{{SL}_2} :=
\Big\{\, b^\beta a^\alpha d^\delta c^\kappa \,\Big|\, \beta, \alpha, \delta,
\kappa \in \N \, ,  \, 0 \in \{\alpha, \delta\} \Big\} \; $  as PBW-like
basis over  $ \Qq \, $.
%
%
%
%
   The definition of the integral form  $ \gerfqsl \, $,  as well as its
proper-\break
 ties, are exactly like those of  $ \gerfqgl \, $,  up to switching
``$ \gergl $''  with  ``$ \gersl $''  and  ``$ GL $''  with  $ SL $''.
%
%
%
%

\vskip2pt

   Another description of  $ \gerfqm $,  $ \gerfqgl $  and  $ \gerfqsl $  is
possible.  Indeed, using a characterization as algebra of matrix coefficients,
$ \fqm $  naturally embeds into  $ \, \uqgl^* \, $.
%
%
 In particular,
there is a perfect (\,=\, non-degenerate) Hopf pairing between  $ \fqm $
and  $ \uqgl $,  which we denote by  $ \; \langle \ , \ \rangle \, \colon
\fqm \times \uqgl \loongrightarrow \Qq \; $  (see e.g.~[No] for details).
Then  $ \; \gerfqm \, = \, \left\{ f \in \fqm \;\Big|\; \big\langle f,
\geruqgl \big\rangle \subseteq \Zqqm \right\} \; $.
%
%
 This leads us to define (cf.~[Ga1])
 \vskip-5pt
  $$  \calfqm \; := \; \left\{ f \in \fqm \;\Big|\; \big\langle f,
\caluqgl \big\rangle \subseteq \Zqqm \right\} \; .  $$
 \vskip-3pt
   The arguments in [Ga1],  {\it mutatis mutandis},  prove also that  $ \; \Q
\otimes_\Z \calfqm \; $  is a  $ \Qqqm $--integral  form of  $ \, \fqm \, $.
Moreover, the analysis therein together with [Ga3], \S\; 7.10, proves that
$ \; \Q \cdot \calfqm \; $  is a  {\sl quantization\/}  of  $ \, U\big( {\gergl_2}^{\!\!*} \big) \, $,  \, i.e.~$ \; {\big( \Q \otimes \calfqm
\big)}_1 \cong U \big( {\gergl_2}^{\!\!*} \big) \; $  as co-Poisson bialgebras
(taking on left-hand side the co-Poisson structure inherited from  $ \; \Q
\otimes \calfqm \, $).

\vskip2pt

   The perfect Hopf pairing between  $ \fqm $  and  $ \uqgl $  uniquely
extends to a similar pairing  $ \; \langle \ , \ \rangle \, \colon \fqgl
\times \uqgl \longrightarrow \Qq \, $,  \, and  $ \; \gerfqgl \, = \left\{
f \in \fqgl \;\Big|\; \big\langle f, \geruqgl \big\rangle \subseteq \Zqqm
\right\} \, $.  This gives the idea (like in [Ga1]) to define
 \vskip-5pt
  $$  \calfqgl \; := \; \left\{\, f \in \fqg \;\Big|\; \big\langle f,
\caluqg \big\rangle \subseteq \Zqqm \right\} \; .  $$
 \vskip-3pt
   \indent   Again,  $ \; \Q \otimes_\Z \calfqgl \; $  is a  $ \Qqqm $--integral
form of  $ \, \fqgl \, $,  \, and that  $ \; \Q \otimes_\Z \calfqgl \; $  is a
{\sl quantization\/}  of  $ U\big({\gergl_2}^{\!\!*}\big) $,  i.e.~$ \; {\big(
\Q \otimes_\Z \calfqgl \big)}_1 \cong U \big( {\gergl_2}^{\!\!*} \big) \; $  as
co-Poisson Hopf algebras.
%
%

\vskip2pt

   The pairing between  $ \fqm $  (or  $ \fqgl $)  and  $ \uqgl $  induces
a perfect pairing between  $ \fqsl $  and  $ \uqsl \, $,  \, giving  $ \;
\gerfqsl \, = \, \Big\{ f \! \in \! \fqsl \;\Big|\, \big\langle f, \geruqsl
\big\rangle \subseteq \Zqqm \Big\} \; $.
%
%
 Then
 \vskip-6pt
  $$  \calfqsl \; := \; \Big\{\, f \in \fqsl \;\Big|\; \big\langle f,
\caluqsl \big\rangle \subseteq \Zqqm \Big\}  $$
 \vskip-1pt
 \noindent   and similar results to those for  $ \calfqgl $  hold for
$ \calfqsl $  too   --- see [Ga1].
%
%

 \vskip2pt

   Finally, let  $ \, \ell \in \N_+ \, $  be odd, and let  $ \, \varepsilon
\, $  be a (formal) primitive  $ \ell $-th  root of 1 as in \S\; 1.3.  Set
$ \; \calfem := {\big(\calfqm\big)}_\varepsilon \, $:  \, then again [Ga1]
and [Ga3] show there is an epimorphism
 \vskip-6pt
  $$  {{\Cal F}r}_{M_2}^{\,\Q} \, \colon \, \Qeps \otimes_{\Zeps}
\calfem \llongtwoheadrightarrow \, \Qeps \otimes_\Z U_\Z \big(
{\gergl_2}^{\!\!*} \big) \, = \, \Qeps \otimes_\Q U \big(
{\gergl_2}^{\!\!*} \big)   \eqno (2.1)  $$
 \vskip-2pt
\noindent
 of bialgebras, which we call  {\sl quantum Frobenius morphism\/}
for  $ \, {\gergl_2}^{\!\!*} \, $.  Similarly, there exist
 \vskip-6pt
  $$  {\calF{}r}_{{GL}_2}^{\,\Q} \, \colon \, \Qeps \otimes_{\Zeps}
\calfegl \llongtwoheadrightarrow \, \Qeps \otimes_\Z \calfunogl \,
\cong \, \Qeps \otimes_\Q U\big({\gergl_2}^{\!\!*}\big)   \eqno (2.2)  $$
 \vskip-4pt
\noindent
  a Hopf algebra epimorphism extending  $ {\calF{}r}_{M_2}^{\,\Q} $,  the
{\sl quantum Frobenius morphism\/}  for  $ \, {\gergl_2}^{\!*} \, $,  \, and
a Hopf algebra epimorphism, uniquely induced by  $ {\calF{}r}_{M_2}^{\,\Z} \! $
or  $ {\calF{}r}_{{GL}_2}^{\,\Z} \, $,
 \vskip-8pt
  $$  {\calF{}r}_{{SL}_2}^{\,\Q} \, \colon \, \Qeps \otimes_{\Zeps}
\calfesl \llongtwoheadrightarrow \, \Qeps \otimes_\Z U_\Z \big(
{\gersl_2}^{\!\!*} \big) \, = \, \Qeps \otimes_\Q U \big(
{\gersl_2}^{\!\!*} \big)   \eqno (2.3)  $$
 \vskip-2pt
\noindent
  where  $ \; \calfesl := {\big(\calfqsl\big)}_\varepsilon \, $,  \, which we
call  {\sl quantum Frobenius morphism\/}  for  $ \, {\gersl_2}^{\!\!*} \, $.

\vskip4pt

   By construction a bialgebra and a Hopf algebra epimorphism  $ \, \fqm
\relbar\joinrel\relbar\joinrel\twoheadrightarrow \fqsl \, $  and  $ \,
\fqgl \relbar\joinrel\relbar\joinrel\twoheadrightarrow \fqsl \, $  exist,
dual to  $ \, \uqsl \longhookrightarrow \uqgl \, $,  \, and similarly
there are epimor\-phisms  $ \, \gerfqm
\relbar\joinrel\relbar\joinrel\twoheadrightarrow \gerfqsl \, $  and
$ \, \gerfqgl \relbar\joinrel\relbar\joinrel\twoheadrightarrow \gerfqsl
\, $  dual to  $ \, \geruqsl \longhookrightarrow \geruqgl \; $.

\vskip11pt

  {\bf 2.3  Dual quantum enveloping algebras.} \, The linear dual
$ {\uqsl}^* $  of  $ \uqsl $  can be seen again (cf.~[Ga1]) as a
quantum group on its own: indeed, we set  $ \, \uqsls := {\uqsl}^* $,
\, a notation used because  $ \uqsls $  stands for the Lie bialgebra
$ {\gersl_n}^{\!\!*} $  just like  $ \uqsl $  stands for  $ \gersl_2 \, $.
Namely,  $ \uqsls $  is a topological Hopf  $ \Qq $--algebra,  with two
integral forms  $ \geruqsls $  and  $ \caluqsls $  which play for
$ \uqsls $  the same r\^{o}le as  $ \geruqsl $  and  $ \caluqsl $
for  $ \uqsl \, $.

\vskip1pt

   The construction goes as follows.  Let  $ \, \H^g_q \, $  be the unital
associative  $ \Qq $--algebra  with generators  $ \, F $,  $ \Lambda_1^{\pm 1} $,
$ \Lambda_2^{\pm 1} $,  $ E $  and relations
 \vskip-17pt
  $$  \displaylines{
   E F = F E  \, ,  \quad \;\;  \Lambda_i^{\,-1} \Lambda_i = 1 = \Lambda_i
\Lambda_i^{\,-1} \, ,  \quad  \Lambda_i^{\pm 1} \Lambda_j^{\pm 1} =
\Lambda_j^{\pm 1} \Lambda_i^{\pm 1} \, ,  \quad  \Lambda_i^{\mp 1}
\Lambda_j^{\pm 1} = \Lambda_j^{\pm 1} \Lambda_i^{\mp 1}
\quad \;  \forall \; i, j  \cr
   \Lambda_1^{\,\pm 1} E = q^{\pm 1} E \Lambda_1^{\,\pm 1} \, ,
\quad  \Lambda_1^{\,\pm 1} F = q^{\pm 1} F \Lambda_1^{\,\pm 1} \, ,
\quad  \Lambda_2^{\,\pm 1} E = q^{\mp 1} E \Lambda_2^{\,\pm 1} \, ,  \quad
\Lambda_2^{\,\pm 1} F = q^{\mp 1} F \Lambda_2^{\,\pm 1} \; .  \cr }  $$
 \vskip-5pt
\noindent
 Let also  $ \H^s_q $  be the
%
%
%
 obtained from  $ \H^g_q $  adding the relation  $ \, \Lambda_1 \Lambda_2
\! = \! 1 \, $.
                                            \par
   The set of PBW-like ordered monomials  $ \; B^g_* := \big\{ E^\eta
\Lambda_1^{\lambda_1} \Lambda_2^{\lambda_2} F^\varphi \,\big|\, \eta,
\lambda_1, \lambda_2, \varphi \in \N \,\big\} \, $  is a  $ \Qq $--basis
for  $ \H^g_q \, $;  similarly  $ \, B^s_* := \big\{ E^\eta \Lambda_1^{\lambda_1}
F^\varphi \,\big|\, \eta, \varphi \! \in \! \N \, , \lambda_1 \in \Z \,\big\} \, $
is a  $ \Qq $--basis  for  $ \H^s_q \; $.
                                            \par
   One defines  $ \uqsls $  as a suitable completion of  $ \H^s_q \, $,
so that  $ \uqsls $  is a topological  $ \Qq $--algebra  topologically
generated by  $ \H^s_q \, $,  and  $ B^s_* $  is a $ \Qq $--basis  of
$ \uqsls $  in topological sense.  Then  $ \uqsls $  is also a topological
Hopf  $ \Qq $--algebra (see [Ga1]).  The same construction makes sense
with  $ \H^g_q $  instead of  $ \H^s_q $  and yields the definition of
$ \uqgls $,  a topological Hopf algebra with  $ B^g_* $  as (topological)
$ \Qq $--basis.  Then by construction  $ \uqsls $  is a quotient of
$ \uqgls \, $,  \, as a topological Hopf algebra, via  $ \; \uqsls_q
\cong \, \uqgls \Big/ ( \Lambda_1 \Lambda_2 - 1) \; $.
                                                \par
   The restricted integral form  $ \geruqsls $  of  $ \uqsls $  is,
by definition, a dense Hopf  $ \Zqqm $--subalgebra  of the subset of
linear functionals in  $ \uqsls $  which are  $ \Zqqm $--valued  onto
$ \caluqsls \, $.  In order to describe it,
set  $ \, L^{\pm 1} := \Lambda_1^{\pm 1} \, $,  and let  $ \,
\gerH^s_q \, $  be the  $ \Zqqm $--subalgebra  of  $ \H^s_q $
generated by all the  $ F^{(m)} $'s,  $ E^{(m)} $'s,  $ L^{\pm 1} $
and  $ \Big(\! {{L \, ; \, c} \atop n} \!\Big) $'s  (for  $ \, m \in
\N \, $,  $ c \in \Z \, $):  \, then  $ \, \frak{B}^s_* := \Big\{\,
E^{(\eta)} \Big(\! {{L \, ; \, 0} \atop l} \!\Big) \, L^{-\text{\it
Ent\,}(l/2)} F^{(\varphi)} \;\Big|\; \eta, l, \varphi \! \in \! \N \,\Big\}
\, $  is a  $ \Zqqm $--basis  of  $ \gerH^s_q \, $,  while  $ \geruqsls $
is the topological closure of  $ \gerH^s_q \, $,  \, and  $ \frak{B}^s_* $
is a topological  $ \Zqqm $--basis  of  $ \geruqsls \, $.
                                                \par
   Similarly, for the integral form  $ \geruqgls $
of  $ \uqgls \, $,  \, take the  $ \Lambda_h $'s  ($ h = 1, 2 $)
instead of  $ L $  and  $ \, \frak{B}^g_* := \Big\{\, E^{(\eta)}
\Big(\! {{\Lambda_1 ; \, 0} \atop \lambda_1} \!\Big) \,
\Lambda_1^{-\text{\it Ent\,}({\lambda_1}/2)} \Big(\! {{\Lambda_2 ;
\, 0} \atop \lambda_2} \!\Big) \, \Lambda_2^{-\text{\it Ent\,}({\lambda_2}
/2)} F^{(\varphi)} \;\Big|\; \eta, \lambda_1, \lambda_2, \varphi \! \in
\! \N \,\Big\} \, $  instead of  $ \frak{B}^s_* \, $.
 By construction  $ \gerH^s_q $  is a quotient algebra of  $ \gerH^g_q $
--- restricting  $ \, \H^g_q \Big/ \! ( \Lambda_1 \Lambda_2 - 1) \cong
\H^s_q \, $  to  $ \gerH^g_q \, $  ---   so  $ \geruqsls $  is a quotient
of  $ \geruqgls \, $,  as a topological Hopf algebra.
 \vskip3pt
   We can describe  $ \gerH^s_q $  rather explicitly: it is the unital
associative  $ \Zqqm $--algebra  with generators  $ \, F^{(m)} $,
$ E^{(m)} $,  $ L^{\pm 1} $,  $ \Big(\! {{L \, ; \, c} \atop m} \!\Big)
\, $   --- for  $ \, m \in \N \, $,  $ c \in \Z $  ---   and relations
 \vskip-10pt
  $$  \displaylines{
   \hbox{\sl relations\/  {\rm (1.1)}  for}  \;\, X = L \; ,
  \hskip15pt  L \, L^{-1} \, = \, 1 \, = \, L^{-1} \, L \; ,
  \hskip15pt  \hbox{\sl relations\/  {\rm (1.2)}  for}  \;\,
X \in \big\{ F \, , E \big\}  \cr
   L^{\pm 1} \, F^{(m)} \, = \; q^{\pm m} \, F^{(m)} \, L^{\pm 1} \; ,
\hskip21pt  E^{(r)} F^{(s)} \, = \; F^{(s)} E^{(r)} \; ,  \hskip21pt
L^{\pm 1} \, E^{(m)} \, = \; q^{\pm m} \, E^{(m)} \, L^{\pm 1}  \cr
%
%
   \textstyle{ \Big(\! {{L \, ; \; c} \atop t} \!\Big) } \, E^{(m)} \, = \, E^{(m)}
\, \textstyle{ \Big(\! {{L \, ; \; c + m} \atop t} \!\Big) } \;\; ,  \hskip41pt   \textstyle{ \Big(\! {{L \, ; \; c} \atop t} \!\Big) } \, F^{(m)} \, = \, F^{(m)}
\, \textstyle{ \Big(\! {{L \, ; \; c + m} \atop t} \!\Big) } \;\; .  \cr }  $$
 \vskip-3pt
%
%
   Similarly,  $ \gerH^g_q $  is the unital associative
$ \Zqqm $--algebra  with generators  $ \, F^{(m)} $,  $ E^{(m)} $,
$ \Lambda_k^{\pm 1} $,  $ \Big(\! {{\Lambda_k \, ; \; c} \atop m}
\!\Big) \, $  (for  $ \, m \in \N \, $,  $ \, c \in \Z \, $,
$ \, k \in \{1,2\} \, $)  and relations
 \vskip-11pt
  $$  \displaylines{
   \Lambda_k \Lambda_k^{-1} \, = \, 1 \, = \, \Lambda_k^{-1} \Lambda_k \; ,
 \hskip27pt  \textstyle{ \Big(\! {{\Lambda_h \, ; \; c} \atop m} \!\Big) \Big(\! {{\Lambda_k \, ; \; s} \atop n} \!\Big) } \, = \, \textstyle{ \Big(\! {{\Lambda_k \, ; \; s} \atop n} \!\Big) \Big(\! {{\Lambda_h \, ; \; c} \atop m} \!\Big) }  \cr
%
%
   \hbox{\sl relations\/  {\rm (1.1)}  for all}
\;\, X \in \big\{\Lambda_1,\Lambda_2\big\} \; ,
\hskip27pt  \hbox{\sl relations\/  {\rm (1.2)}  for all}
\;\, X \in \big\{ F \, , E \big\}  \cr
   E^{(r)} F^{(s)} \, = \; F^{(s)} E^{(r)} \; ,  \hskip17pt
\Lambda_k^{\pm 1} \, Y^{(m)} = \, q^{\pm (\delta_{k,1} -
\delta_{k,2}) m} \, Y^{(m)} \Lambda_k^{\pm 1}  \hskip17pt
\forall \;\, Y \in \big\{F \, , E \big\}  \cr
%
%
   \textstyle{ \Big(\! {{\Lambda_k \, ; \; c} \atop t} \!\Big) } \, Y^{(m)}
\, = \, Y^{(m)} \, \textstyle{ \Big(\! {{\Lambda_k \, ; \; c + (\delta_{k,1}
- \delta_{k,2}) \, m} \atop t} \!\Big) }   \hskip23pt  \forall \;\, Y \in
\big\{ F \, , E \big\} }  $$
 \vskip-1pt
%
%
 In this paper we do not need the Hopf structure of  $ \geruqsls $
and  $ \geruqsls $  (cf.~[Ga1]).
 \vskip2pt
   $ \geruqsls $  is a  {\sl quantization\/}  of  $ U_\Z\big(
{\gersl_2}^{\!\!*} \big) $,  for  $ \, \geruunosls \! := \! {\big(
\geruqsls \big)}_1 \! \cong \! U_\Z \big({\gersl_2}^{\!\!*}\big) \, $
as co-Poisson Hopf algebras, with on left-hand side the co-Poisson
structure inherited from  $ \geruqsls \, $.  In terms of generators
(notation of \S\; 1.2) this reads  $ \, F^{(m)}{\big|}_{q=1} \! \cong
\text{f}^{\,(m)} \, $,  $ \, \Big(\! {{L \, ; \, 0} \atop m} \!\Big)
{\Big|}_{q=1} \! \cong \Big(\! {\text{h} \atop m} \!\Big) \, $,  $ \,
L^{\pm 1}{\big|}_{q=1} \! \cong 1 \, $,  $ \, E^{(m)}{\big|}_{q=1} \!
\cong \text{e}^{(m)} \, $  for  $ \, m \in \N \, $.   Similarly  $ \,
\geruunogls \! \cong U_\Z \big({\gergl_2}^{\!\!*}\big) \, $,  with
$ \, F^{(m)}{\big|}_{q=1} \! \cong \text{f}^{\,(m)} $,  $ \, \Big(\!
{{\Lambda_k \, ; \, 0} \atop m} \!\Big){\Big|}_{q=1} \! \cong \Big(\!
{\text{g}_k \atop m} \!\Big) \, $,  $ \, \Lambda_k^{\pm 1}{\big|}_{q=1}
\! \cong 1 \, $,  $ \, E^{(m)}{\big|}_{q=1} \! \cong \text{e}^{(m)} \, $
for  $ \, m \in \N \, $  and  $ \, k \in \{1,2\} \, $.
%
%
                                          \par
   Finally, let  $ \ell $  and  $ \varepsilon $  be as in \S\; 1.3.  Set
$ \, \geruesls := {\big(\geruqsls\big)}_\varepsilon \, $  and  $ \,
\gerH^s_\varepsilon := {\big(\gerH^s\big)}_\varepsilon \, $.  Then
(cf.~[Ga1], \S\; 7.7) the embedding  $ \, \gerH^s_\varepsilon
\longhookrightarrow \geruesls \, $  is an isomorphism, thus  $ \,
\geruesls = \gerH^s_\varepsilon \, $.  Similarly (with like notation)
$ \, \geruegls = \gerH^g_\varepsilon \, $.  Also, there are Hopf
algebra epimorphisms
 \vskip-11pt
  $$  \eqalignno{
   \gerFr_{{\gersl_2}^{\!\!*}}^{\,\Z} \, \colon \, \geruesl =
\gerH^s_\varepsilon \, \llongtwoheadrightarrow \, \Zeps \otimes_\Z
\gerH^s_1  &  = \Zeps \otimes_\Z \geruunosls \, \cong \, \Zeps
\otimes_\Z U_\Z\big({\gersl_2}^{\!\!*}\big)  \qquad   &  (2.4)  \cr
   \gerFr_{{\gergl_2}^{\!\!*}}^{\,\Z} \, \colon \, \geruegl =
\gerH^g_\varepsilon \, \llongtwoheadrightarrow \, \Zeps \otimes_\Z
\gerH^g_1  &  = \Zeps \otimes_\Z \geruunogls \, \cong \, \Zeps
\otimes_\Z U_\Z\big({\gergl_2}^{\!\!*}\big)  \qquad   &  (2.5)  \cr }  $$
%
%
 \eject
\noindent
 defined by  $ \, X^{(s)}\big\vert_{q=\varepsilon} \!\!\! \mapsto
\text{x}^{\, (s / \ell)} \, $,  $ \, \Big(\! {{Y \, ; \; 0} \atop s}
\!\Big) \! \Big\vert_{q=\varepsilon} \!\!\! \mapsto \Big(\! {{\text{k}}
\atop {s / \ell}} \!\Big) \, $  if  $ \, \ell \Big\vert s \, $,  $ \,
X^{(s)} \big\vert_{q=\varepsilon} \!\!\! \mapsto 0 \, $,  $ \, \Big(\!
{{Y \, ; \; 0} \atop s} \!\Big) \! \Big\vert_{q=\varepsilon} \!\!\!
\mapsto 0 \, $  if  $ \, \ell  \hbox{$ \not\!\big\vert $} s \, $,
and  $ \, Y^{\pm 1} \Big\vert_{q=1} \!\!\! \mapsto 1 \, $,  with  $ \,
\big(X,\text{x}\big) \in \big\{\! \big(F,\text{f}\,\big), \big( E,
\text{e} \big) \!\big\} \, $,  \, and  $ \, \big(Y,\text{k}\big) =
\big(L,\text{h}\big) \, $  in the  $ \gersl_n $  case,  $ \, \big(Y,
\text{k}\big) \in \big\{\!\big(\Lambda_i,\text{g}_i\big)\!\big\}_{i=1,2}
\, $  for  $ \gergl_n \, $.  These are  {\it quantum Frobenius morphism\/}
for  $ {\gersl_n}^{\!\!*} $  and  $ {\gergl_n}^{\!\!*} $.
                                           \par
   The above epimorphism  $ \; \pi_q \, \colon \, \geruqgls
\twoheadrightarrow \geruqsls \, $  of topological Hopf
$ \Zqqm $--algebras  is compatible with the
quantum Frobenius morphisms,
%
%
in the obvious sense.
                                           \par
   To finish with, the natural evaluation pairing  $ \; \langle \,\ ,
\ \rangle \, \colon \uqgs \times \uqg \longrightarrow \Qq \; $  (for
$ \, \gerg \in \{ \gersl_2, \gergl_2 \} \, $)  is uniquely determined
by its values on PBW bases: we have them via
 \vskip-11pt
  $$  {\textstyle \Big\langle
    E^{(\eta)} \, \Big(\! {{L \, ; \, 0} \atop l} \!\Big)
\, L^{-\text{\it Ent}\,(l/2)} \, F^{(\varphi)} \; , \;
   \fbar^{\,f} } \, {\textstyle K^\kappa \, \ebar^{\,e} \Big\rangle
\; = \; {(-1)}^\eta \, \delta_{e,\varphi} \, \delta_{f,\eta} \,
\Big(\! {\kappa \atop l} \!\Big)_{\!q} \, q^{- \kappa
\text{\it Ent}\,(l/2)} }   \eqno (2.6)  $$
 \vskip-17pt
  $$  \hbox{ $ \eqalign{
         {\textstyle \Big\langle
    E^{(\eta)} \, \Big(\! {{\Lambda_1 \, ; \, 0} \atop \lambda_1} \!\Big)
\, \Lambda_1^{-\text{\it Ent\,}(\lambda_1 / 2)} \Big(\! {{\Lambda_2 \, ;
\, 0} \atop \lambda_2} \!\Big) \, \Lambda_2^{-\text{\it Ent\,}
(\lambda_2 / 2)} \, F^{(\varphi)} \; , \;
   \fbar^{\,f} }  &  {\textstyle G_1^{\!\phantom{|}\gamma_1} \,
{\textstyle G_2^{\!\phantom{|}\gamma_2} \, \ebar^{\,e} \Big\rangle }
\; = }  \cr
   {\textstyle = \; {(-1)}^\eta \, \delta_{e,\varphi} \, \delta_{f,\eta}
\, \Big(\! {\gamma_1 \atop \lambda_1} \!\Big)_{\!q} \Big(\! {\gamma_2
\atop \lambda_2} \!\Big)_{\!q} \, }  &  {\textstyle q^{- \gamma_1
\text{\it Ent}\,(\lambda_1 / 2) - \gamma_2 \text{\it Ent}\,
(\lambda_2 / 2)} }  \cr } $ }   \eqno (2.7)  $$

\vskip9pt

  {\bf 2.4 Embedding quantum function algebras into quantum enveloping
algebras.} \, Let  $ \, G \in \big\{ {SL}_2, {GL}_2 \big\} \, $  and
$ \, \gerg := \text{\it Lie}\,(G) \, $.  By definition  $ \fqg $  embeds
into  $ \, \uqgs := {\uqg}^* \, $,  \, via a monomorphism  $ \; \xi \,
\colon \, \fqg \longhookrightarrow \uqgs \; $  of
        topological Hopf  $ \Qq $--algebras.\break
\noindent
 Moreover  $ \; \calfqg = \xi^{-1}\big(\geruqgs\big) \, $,  \; so  $ \xi $  restricts to a monomorphism  $ \; \widehat{\xi} \, \colon \, \calfqg
\longhookrightarrow \geruqgs \; $  too, and similarly  $ \, \widetilde{\xi}
\, \colon \, \gerfqg \longhookrightarrow \caluqgs \, $.  These verify
$ \, \xi\big(\fqg\big) \subseteq \H_q^x \, $  and  $ \,
\widehat{\xi}\big(\calfqg\big) \subseteq \gerH_q^{\,x} \, $  (with  $ \,
x \in \{s,g\} \, $,  according to the type of  $ G \, $)  so  $ \, \fqg
= \xi^{-1}(\H_x) \, $  and  $ \, \calfqg \, = \, \widehat{\xi}^{-1}(\gerH_x)
\, $.  Furthermore,  $ \widehat{\xi} $  is compatible with specializations
and quantum Frobenius morphisms, that is  $ \; \Big( \text{id}_{\Qeps}
\otimes_{\Z} \widehat{\xi} {\,\big|}_{q=1} \Big) \circ \calFr_G^{\,\Q}
= \Big( \text{id}_{\Qeps} \otimes_{\Zeps} \gerFr_{\gerg^*}^{\,\Z} \Big)
\circ \Big( \text{id}_{\Qeps} \otimes_{\Zeps} \widehat{\,\xi}
{\big|}_{q=\varepsilon} \Big) \; $.  As  $ \fqm $  embeds into
$ \fqgl $,  restricting  $ \, \xi \, \colon \fqgl \lhook\joinrel\rightarrow
\uqgls \, $  yields an embedding  $ \, \xi \, \colon \fqm
\lhook\joinrel\rightarrow \uqgls \, $,  \; and similarly we have
an embedding  $ \; \widehat{\xi} \, \colon \, \calfqm
\longhookrightarrow \geruqgls \, $,  \; which factors through
$ \calfqgl $  (and similarly  $ \; \widetilde{\xi} \, \colon \,
\gerfqm \longhookrightarrow \caluqgls \, $,  \, which factors
through  $ \gerfqgl \, $).
                                            \par
   In [Ga1], Appendix, embeddings  $ \xi $  and  $ \widetilde{\xi} $
as above are described for  $ {SL}_2 $,  namely given by  $ \; \xi \,
\colon \, a \mapsto L - \fbar L^{-1} \ebar \, $,  $ \, b \mapsto -
\fbar L^{-1} \, $,  $ \, c \mapsto + L^{-1} \ebar \, $,  $ \, d
\mapsto L^{-1} \; $.  Similarly (and exploiting the analysis
in [Ga2], \S\S\; 5.2/4), we find an analogous  $ \xi $  for
$ {GL}_2 $,  namely  $ \; \xi \, \colon \, \fqgl \longhookrightarrow
\uqgls \, $,  $ \, a \mapsto \Lambda_1 - \fbar \Lambda_2
\ebar \, $,  $ \, b \mapsto - \fbar \Lambda_2 \, $,
$ \, c \mapsto + \Lambda_2 \ebar \, $,  $ \, d \mapsto
\Lambda_2 \, $,  $ \, {D_q}^{\!-1} \mapsto \big( \Lambda_1
\Lambda_2 \big)^{-1} \; $.  The same formul\ae{}  describe $ \;
\widetilde{\xi} \, \colon \, \gerfqgl \longhookrightarrow
\caluqgls \, $.  Discarding  $ \, {D_q}^{\!-1} \! \mapsto \big(
\Lambda_1 \Lambda_2 \big)^{-1} \, $  they describe also the
embedding  $ \; \xi \, \colon \, \fqm \longhookrightarrow \uqgls
\, $,  $ \, a \mapsto \Lambda_1 - \fbar \Lambda_2 \ebar
\, $,  $ \, b \mapsto - \fbar \Lambda_2 \, $,  $ \, c \mapsto
+ \Lambda_2 \ebar \, $,  $ \, d \mapsto \Lambda_2 \, $,
obtained restricting  $ \; \xi \, \colon \, \fqgl \longhookrightarrow
\uqgls \, $,  \; and its restriction  $ \; \widetilde{\xi} \, \colon
\, \gerfqm \longhookrightarrow \caluqgls \, $.
                                            \par
   Finally, the various embeddings  $ \xi $  and their
restrictions to integral forms also are compatible   ---
in the obvious sense ---   with the epimorphisms  $ \,
\uqgls \relbar\joinrel\relbar\joinrel\twoheadrightarrow
\uqsls \, $  and  $ \, \fqgl
\relbar\joinrel\relbar\joinrel\twoheadrightarrow
\fqsl \, $  or  $ \, \fqm
\relbar\joinrel\relbar\joinrel\twoheadrightarrow
\fqsl \, $  and their restrictions to integral forms.
%
%

 \vskip1,1truecm

\centerline {\bf \S\; 3 \ The structure of  $ \, \calfqm \, $,  $ \,
\calfqgl \, $  and  $ \, \calfqsl \, $,}
\centerline {\bf their specializations and quantum Frobenius epimorphisms. }

\vskip9pt

   We need some more notation.  First set  $ \; \b :=
\big( q - q^{-1} \big)^{-1} b \; $  and  $ \; \c := \big( q
- q^{-1} \big)^{-1} c \, $.  Then for all  $ \, n \in \N \, $,
\, we set  $ \; \b^{(n)} \! := \b^n \! \big/ [n]_q! \, $,  $ \;
\c^{(n)} \! := \c^n \! \big/ [n]_q! \; $   --- like in \S\; 1.2.

\vskip9pt

\proclaim{Theorem 3.1}
 \vskip3pt
   (a) \,  $ \calfqm $  is a free  $ \, \Zqqm $--module,  with basis the
set of ordered monomials
 \vskip-5pt
  $$  \Cal{B}_{M_2} \, = \, \textstyle{ \left\{\, \b^{(\beta)} \, \Big(\!
{{a \, ; \, 0} \atop \alpha} \!\Big) \, \Big(\! {{d \, ; \, 0} \atop \delta}
\!\Big) \, \c^{(\kappa)} \;\Big|\; \alpha, \beta, \kappa, \delta \in \N
\,\right\} }  $$
 \vskip-1pt
\noindent
(a PBW-like basis).  Similarly, any other set obtained from  $ \Cal{B}_{M_2} $
via permutations of factors (of the monomials in  $ \Cal{B}_{M_2} $)
is a  $ \, \Zqqm $--basis  of  $ \, \calfqm $  as well.
 \vskip3pt
   (b) \,  $ \calfqgl $  is the  $ \, \Zqqm $--span  of
the set of ordered monomials
 \vskip-3pt
  $$  \Cal{S}_{{GL}_2} \, = \, \textstyle{ \left\{\, \b^{(\beta)} \, \Big(\!
{{a \, ; \, 0} \atop \alpha} \!\Big) \Big(\! {{d \, ; \, 0} \atop \delta}
\!\Big) \, \c^{(\kappa)} \, {D_q}^{\!-\nu} \;\Big|\; \alpha, \beta, \kappa,
\delta, \nu \in \N \,\right\} } \quad .  $$
 \vskip-1pt
\noindent
 Similarly, any other set obtained from  $ \Cal{S}_{{GL}_2} $  via
permutations of factors (of the monomials in  $ \Cal{S}_{{GL}_2} $)
is a  $ \, \Zqqm $--spanning set for  $ \, \calfqgl \, $  as well.
Moreover, if  $ \, f \in \calfqgl \, $  then  $ f $  can be expanded into
a  $ \Zqqm $--linear  combination of elements of  $ \Cal{S}_{{GL}_2} $
which all bear the same exponent  $ \nu \, $;  \, similarly for the other
spanning sets mentioned above.
 \vskip3pt
   (c) \,  $ \calfqsl $  is the  $ \, \Zqqm $--span  of the
set of ordered monomials
 \vskip-3pt
  $$  \Cal{S}_{{SL}_2} \, = \, \textstyle{ \left\{\, \b^{(\beta)} \, \Big(\!
{{a \, ; \, 0} \atop \alpha} \!\Big) \Big(\! {{d \, ; \, 0} \atop \delta}
\!\Big) \, \c^{(\kappa)} \;\Big|\; \alpha, \beta, \kappa, \delta \in \N
\,\right\} } \quad .  $$
 \vskip-1pt
\noindent
 Similarly, any other set obtained from  $ \Cal{S}_{{SL}_2} $  via
permutations of factors (of the monomials in  $ \Cal{S}_{{SL}_2} $)
is a  $ \, \Zqqm $--spanning set for  $ \, \calfqsl \, $  as well.
\endproclaim

\demo{Proof}  {\it (a)} \,  For all  $ \, \alpha, \beta, \kappa, \delta
\in \N \, $,  \, let  $ \; \Cal{M}_{\alpha,\beta,\kappa,\delta} \, ;= \,
\b^{(\beta)} \, {{a \, ; \, 0} \choose \alpha} \, {{d \, ; \, 0} \choose
\delta} \, \c^{(\kappa)} \, $.  Due to the formul\ae{}  for  $ \; \xi \,
\colon \, \fqm \longhookrightarrow \uqgls \, $  in \S\; 2.4 and to  Lemma
4.2{\it (b)\/}  later on, one has
 \vskip-15pt
  $$  \displaylines{
 \textstyle
   \xi \big( \Cal{M}_{\alpha,\beta,\kappa,\delta} \big) \, = \,
 \big(- \overline{F} \Lambda_2 \big)^{(\beta)} \cdot
\Big(\! {{\Lambda_1 - \overline{F} \Lambda_2 \overline{E} \, ; \, 0}
\atop \alpha} \Big) \cdot \Big(\! {{\Lambda_2 \, ; \, 0} \atop \delta}
\Big) \cdot \big(+ \Lambda_2 \overline{E} \big)^{(\kappa)} \, = \,
\sum_{r=0}^\alpha \, {(-1)}^{\beta + r} \times   \hfill  \cr
 \textstyle
   \; \times \, q^{{\kappa \choose 2} - {\beta \choose 2} - {r \choose 2}
- r (\alpha+1)} \big( q - q^{-1} \big)^r \, {[r]}_q! \, {\left[ {{\beta
+ r} \atop \beta} \right]}_{\!q} \, {\Big[ {{\kappa + r} \atop \kappa}
\Big]}_{\!q} F^{(\beta + r)} \left\{ {\Lambda_1 \, ; \, 0} \atop
{\alpha \; , \; r} \right\} \Big(\! {{\Lambda_2 \, ; \, 0} \atop
\delta} \Big) \, {\Lambda_2}^{\beta + r + \kappa} \, E^{(\kappa + r)}
\cr }  $$
 \vskip-6pt
\noindent
 so that  $ \; \xi \big( \Cal{M}_{\alpha,\beta,\kappa,\delta} \big)
\in \gerH^g_q \, $,  \; which implies, thanks to  $ \, \calfqg \,
= \, \widehat{\xi}^{-1}(\gerH_x) \, $  (see \S\; 2.4),
 \vskip-5pt
  $$  \Cal{M}_{\alpha,\beta,\kappa,\delta} \, \in \,
\widehat{\xi}^{-1}(\gerH^g_q) \, {\textstyle \bigcap} \,
\fqm \, = \, \calfqm \; .   \eqno (3.1)  $$
 \vskip-1pt
\noindent
 This proves that  $ \, \Cal{B}_{M_2} \subseteq \calfqm \, $,  \,
hence the  $ \Zqqm $--span  of  $ \Cal{B}_{M_2} $  is contained
in  $ \calfqm \, $.
                                             \par
   Now pick  $ \, f \in \calfqm \, $.  Clearly  $ \Cal{B}_{M_2} $
is a  $ \Qq $--basis  of  $ \fqm $,  hence there is a unique
expansion  $ \; f = \sum_{\alpha, \beta, \kappa, \delta \in \N}
\chi_{\alpha,\beta,\kappa,\delta} \, \Cal{M}_{\alpha, \beta,
\kappa, \delta} \; $  with all coefficients  $ \, \chi_{\alpha,
\beta, \kappa, \delta} \in \Qq \, $;  \, we must show that
these belong to  $ \Zqqm $.  Let  $ \, \beta_0 \, $  and
$ \kappa_0 $  in  $ \N $  be the least indices such that  $ \,
\chi_{\alpha,\beta_0,\kappa_0,\delta} \not= 0 \, $  for some  $ \,
\alpha, \delta \in \N \, $.  The previous description of  $ \xi
\big( \Cal{M}_{\alpha,\beta,\gamma,\delta} \big) $  and (2.7) yield
 \vskip-3pt
  $$  \hbox{ $ \eqalign{
 \textstyle
   \left\langle \Cal{M}_{\alpha,\beta,\kappa,\delta} \, , \,
\overline{E}^{\,\eta} {G_1}^{\!\gamma_1} {G_2}^{\!
\gamma_2} \, \overline{F}^{\,\varphi} \right\rangle \,  &  = \, 0
\qquad \hskip15pt  \hfill   \text{if}  \quad  \eta < \beta
\text{\ \ or \ }  \varphi < \kappa  \cr
 \textstyle
   \left\langle \Cal{M}_{\alpha,\beta,\kappa,\delta} \, , \,
\overline{E}^{\,\beta} {G_1}^{\!\gamma_1} {G_2}^{\!
\gamma_2} \, \overline{F}^{\,\kappa} \right\rangle \;  &
 \textstyle
= \; {(-1)}^\beta \, q^{{\kappa \choose 2} - {\beta \choose 2}
+ (\beta + \kappa) \gamma_2} \, \Big(\! {\gamma_1 \atop \alpha}
\!\Big)_{\!q} \, \Big(\! {\gamma_2 \atop \delta} \!\Big)_{\!q}
\cr } $ }  $$
 \vskip-2pt
\noindent
  This gives  $ \; \left\langle f \, , \, \overline{E}^{\,\beta_0}
{G_1}^{\!\gamma_1} {G_2}^{\!\gamma_2} \, \overline{F}^{\,\kappa_0}
\right\rangle \, = \,
%
%
 q^{{\kappa_0 \choose 2} - {\beta_0 \choose 2} + (\beta_0 +
\kappa_0) \gamma_2} \, \sum_{\alpha,\delta} \, \chi_{\alpha,\beta_0,
\kappa_0,\delta} \, {\gamma_1 \choose \alpha}_{\!q} \, {\gamma_2
\choose \delta}_{\!q} \; $.
By assumption the last term belongs to  $ \Zqqm $,  whence also
 \vskip-13pt
  $$  y_{\gamma_1, \gamma_2} \; := \; \textstyle \sum\limits_{\alpha,\delta}
\; \chi_{\alpha,\beta_0, \kappa_0,\delta} \; \Big(\! {\gamma_1 \atop \alpha}
\!\Big)_{\!q} \, \Big(\! {\gamma_2 \atop \delta} \!\Big)_{\!q} \, \in \Zqqm
\hskip15pt  \forall \;\; \gamma_1, \gamma_2 \in \Z, \, \alpha, \delta \in
\N \; .   \eqno (3.2)  $$
 \vskip-5pt
\noindent
 Using as indices the pairs  $ \, (\alpha,\delta) \in \N^2 \, $  and
$ \, (\gamma_1,\gamma_2) \in \N^2 \, $,  \, the set of identities (3.2)
can be read as a change of variables from  $ \, \big\{ \chi_{\alpha,
\beta_0, \kappa_0, \delta} \big\}_{(\alpha,\delta) \in \N^2} \, $  to
$ \, \big\{ y_{\gamma_1,\gamma_2} \big\}_{(\gamma_1,\gamma_2) \in \N^2}
\, $;  \, fixing in  $ \N^2 $  any total order  $ \, \preceq \, $  such
that  $ \, (m,n) \preceq (m',n') \, $  if  $ \, m \leq m \, $  or  $ \,
m' \leq n' \, $,  \, the (infinite size) matrix ruling this change of
variables,  i.e.~$ \, {\left( \! {\gamma_1 \choose \alpha}_{\!q} \,
{\gamma_2 \choose \delta}_{\!q} \right)}_{(\gamma_1,\gamma_2), (\alpha,
\delta) \in \N^2} \, $,  \, has entries in  $ \Zqqm $  and is lower
triangular unipotent.  Thus it is invertible and its inverse also is
lower triangular unipotent with entries  in  $ \Zqqm $,  so  $ \,
\chi_{\alpha,\beta_0,\kappa_0,\delta} \in \Zqqm \, $  for all
$ \, \alpha, \delta \in \N \, $.
                                      \par
   The previous analysis gives  $ \, f' := f - \sum_{\alpha,\delta}
\chi_{\alpha, \beta_0,\kappa_0,\delta} \, \Cal{M}_{\alpha,\beta_0,
\kappa_0,\delta} \in \calfqm \, $; moreover, by construction the
expansion of  $ f' $  as a  $ \Qq $--linear  combination of elements
of  $ \Cal{B}_{M_2} $  has less non-trivial summands than  $ f \, $:
\, then we can apply the same argument, and iterate till we find
that all coefficients  $ \chi_{\alpha,\beta,\kappa,\delta} $  in the
original expansion of  $ f $  do belong to  $ \Zqqm \, $.
                                      \par
   Finally, the last observation about other bases is clear.
 \vskip3pt
   {\it (b)} \,  By claim  {\it (a)},  every monomial of type
$ \; \b^{(\beta)} \, \Big(\! {{a \, ; \, 0} \atop \alpha} \!\Big)
\Big(\! {{d \, ; \, 0} \atop \delta} \!\Big) \, \c^{(\kappa)} \; $
is  $ \Zqqm $--integer-valued on  $ \caluqgl \, $;
%
%
 on the other hand, the same is true for  $ \, {D_q}^{\!-\nu} \, $
($ \, \forall \; \nu \in \N \, $),  \, because  $ \; {D_q}^{\!-\nu} \in
\gerfqgl \, $  and  $ \, \gerfqgl \subseteq \calfqgl \, $  since  $ \,
\geruqgl \supseteq \caluqgl \, $.  Finally,
 \vskip-11pt
  $$  \displaylines{ \textstyle {
   \left\langle \b^{(\beta)} \Big(\! {{a \, ; \, 0} \atop \alpha} \!\Big)
\! \Big(\! {{d \, ; \, 0} \atop \delta} \!\Big) \, \c^{(\kappa)} \,
{D_q}^{\!\! -\nu} , \, \caluqgl \! \right\rangle = \left\langle \b^{(\beta)}
\Big(\! {{a \, ; \, 0} \atop \alpha} \!\Big) \! \Big(\! {{d \, ; \, 0}
\atop \delta} \!\Big) \, \c^{(\kappa)} \otimes {D_q}^{\!\! -\nu} ,
\, \Delta\big(\caluqgl\big) \! \right\rangle \subseteq }   \hfill  \cr
   \hfill   \textstyle { \subseteq \; \left\langle \b^{(\beta)} \,
\Big(\! {{a \, ; \, 0} \atop \alpha} \!\Big) \Big(\! {{d \, ; \, 0} \atop \delta} \!\Big) \, \c^{(\kappa)} , \, \caluqgl \right\rangle \cdot \left\langle {D_q}^{\!\! -\nu}, \, \caluqgl \right\rangle \; \subseteq \; \Zqqm } }  $$
 \vskip-7pt
\noindent
 so  $ \, \b^{(\beta)} \Big(\! {{a \, ; \, 0} \atop \alpha} \!\Big)
\Big(\! {{d \, ; \, 0} \atop \delta} \!\Big) \c^{(\kappa)} {D_q}^{\!\! -\nu} \!
\in \! \calfqgl \, $,  and the  $ \Zqqm $--span  of  $ \Cal{S}_{{GL}_2} $  sits
inside  $ \calfqgl \, $.
%
%
                                       \par
   Conversely, let  $ \, f \in \calfqgl \, $.  Then there exists  $ \,
N \in \N \, $  such that  $ \, f \, {D_q}^{\!N} \in \fqm \, $.  In
addition,  $ \; \left\langle f \, {D_q}^{\!N}, \, \caluqgl \right\rangle
\, = \, \left\langle f \otimes {D_q}^{\!N}, \, \Delta\big(\caluqgl\big)
\right\rangle \, \subseteq \left\langle f, \, \caluqgl \right\rangle
\cdot \left\langle {D_q}^{\!N}, \, \caluqgl \right\rangle \, \subseteq
\, \Zqqm \; $  because  $ \, f, {D_q}^{\!N} \in \calfqgl \, $.  Thus
$ \; f \, {D_q}^{\!N} \in \calfqm \, $;  \; then, by claim  {\it (a)\/},
$ \, f \, {D_q}^{\!N} \, $  belongs to the  $ \Zqqm $--span  of  $ \,
\Cal{B}_{M_2} \, $,  \, whence the claim follows at once.
                                      \par
   Finally, the last observation about other spanning sets
is self-evident.
 \vskip3pt
   {\it (c)} \,  The projection epimorphism  $ \; \fqm \,{\buildrel \pi
\over {\llongtwoheadrightarrow}}\, \fqm \Big/ \big(D_q\big) \, \cong
\, \fqsl \; $  (given by restriction from  $ \uqgl $  to  $ \uqsl \, $)
maps  $ \Cal{B}_{M_2} $  onto  $ \Cal{S}_{{SL}_2} \, $:  \, it follows
directly from definitions that  $ \, \pi\big(\calfqm\big) \subseteq
\calfqsl \, $,  \, hence in particular (thanks to claim  {\it (a)\/})
the  $ \Zqqm $--span  of  $ \Cal{S}_{{SL}_2} $  is contained in
$ \calfqsl \, $.  Conversely, let  $ \, f \in \calfqsl \, $.
Since  $ B_{{SL}_2} $  in \S\; 2.2 is a  $ \Qq $--basis  of
$ \fqsl $  it follows that any  $ \, f \in \calfqsl \, $
has a unique expansion
 \vskip-5pt
  $$  \hbox{ $ \eqalign{
   f \;\,  &  = \, {\textstyle \sum_{\alpha, \beta, \kappa \in \N}} \,
\chi^{\,\alpha}_{\beta,\kappa} \cdot \b^{(\beta)} a^\alpha \c^{(\kappa)}
\, + {\textstyle \sum_{\beta, \kappa, \delta \in \N}} \,
\eta^{\,\delta}_{\beta,\kappa} \cdot \b^{(\beta)} d^{\,\delta}
\c^{(\kappa)} \, =  \cr
         &  = \, {\textstyle \sum_{\beta, \kappa \in \N}} \,
\b^{(\beta)} \left( \, {\textstyle \sum_{\alpha \in \N}} \,
\chi^{\,\alpha}_{\beta,\kappa} \, a^\alpha \, + \, \varphi \, + \,
{\textstyle \sum_{\delta \in \N}} \, \eta^{\,\delta}_{\beta,
\kappa} \, d^{\,\delta} \right) \c^{(\kappa)}  \cr } $ }   \eqno (3.3)  $$
 \vskip-5pt
\noindent
 for some  $ \, \chi^{\,\alpha}_{\beta,\kappa}, \eta^{\,\delta}_{\beta,
\kappa}, \varphi \in \Qq \, $.  Since  $ \, D_q := a \, d \, - q \,
b \, c = 1 \, $  in  $ \calfqsl $,  we can rewrite (3.3) as  $ \; f
\, = \, {\textstyle \sum\limits_{\beta, \kappa \in \N}} \, \b^{(\beta)}
\left( \, {\textstyle \sum\limits_{\alpha \in \N}} \, \chi^{\,
\alpha}_{\beta,\kappa} \, a^\alpha \cdot {D_q}^{\!\mu} \, + \,
\varphi \cdot {D_q}^{\!\mu} \, + \, {\textstyle \sum\limits_{\delta
\in \N}} \, \eta^{\,\delta}_{\beta,\kappa} \, d^{\,\delta} \cdot
{D_q}^{\! \mu - \delta} \right) \c^{(\kappa)} \, $,  \; where
$ \; \mu := \max \big\{ \delta \in \N \,\big\vert\, \eta^{\,
\delta}_{\beta,\kappa} \not= 0 \big\} \, $.  Now consider the
element of  $ \fqm $
 \vskip-13pt
  $$  f' \, = \, {\textstyle \sum_{\beta, \kappa \in \N}} \, \b^{(\beta)}
\left( \, {\textstyle \sum_{\alpha \in \N}} \, \chi^{\,\alpha}_{\beta,
\kappa} \, a^\alpha \cdot {D_q}^{\!\mu} \, + \, \varphi \cdot {D_q}^{\!\mu}
\, + \, {\textstyle \sum_{\delta \in \N}} \, \eta^{\,\delta}_{\beta,\kappa}
\, d^{\,\delta} \cdot {D_q}^{\! \mu - \delta} \right) \c^{(\kappa)}  $$
%
%
 \vskip-5pt
\noindent
 By construction,  $ \, \pi(f') = f \, $;  \; moreover, a straightforward
check shows that
 \vskip-9pt
  $$  \left\langle f' \, , \, \overline{E}^{\,\eta} K^\kappa \,
{G_2}^{\!\gamma} \, \overline{F}^{\,\varphi} \right\rangle \, =
\, q^{\gamma \mu} \cdot \left\langle f \, , \, \overline{E}^{\,\eta}
K^\kappa \, \overline{F}^{\,\varphi} \right\rangle \, \in \, \Zqqm  $$
 \vskip-3pt
\noindent
 (for all  $ \, \eta, \varphi \in \N $,  $ \kappa, \gamma \in \Z \, $)
because  $ \, f \in \calfqsl \, $  by hypothesis.  Since the monomials
$ \, \overline{E}^{\,\eta} K^\kappa \, {G_2}^{\!\gamma} \, \overline{F}^{\,
\varphi} \, $  form a  $ \Zqqm $--basis  of  $ \caluqgl $  (by \S\; 2.1,
as  $ \, K := G_1 {G_2}^{\!-1} \, $),  it follows that  $ \, f' \in
\calfqm \, $.  Then claim  {\it (a)\/}  and the fact that  $ \,
\pi(f') = f \, $  imply that  $ f $  lies in the  $ \Zqqm $--span
of  $ \Cal{S}_{{SL}_2} $.  Finally, the last observation
about other bases is clear.   \qed
\enddemo

\vskip11pt

\proclaim{Theorem 3.2}
 \vskip2pt
   (a) \;  $ \calfqm $  is the unital associative  $ \Zqqm $--algebra
with generators
 \vskip-5pt
  $$  {} \hskip25pt   {{a \, ; r} \choose n} \; ,  \quad  \b^{(n)} \; ,
\quad  \c^{(n)} \; ,  \quad  {{d \, ; s} \choose m}   \eqno \forall \;\;
m, n \in \N , \, r, s \in \Z  \hskip25pt {}  $$
 \vskip-3pt
\noindent
 and relations
 \vskip-19pt
  $$  \displaylines{
   \b^{(r)} \c^{(s)} \, = \, \c^{(s)} \b^{(r)}  \cr
   \hbox{\sl relations\/  {\rm (1.1)}  for \ }
X \in \big\{ a , d \big\} \; ,  \qquad  \hbox{\sl relations\/
{\rm (1.2)}  for \ }  X \in \{b,c\}  \cr
%
%
   \left[ {{a \, ; r} \choose n}, {{d \, ; s} \choose m} \right]
\, = \; {\textstyle \sum\limits_{j=1}^{n \wedge m}} \, q^{j \,
((r+s)-(n+m)) + {j \choose 2}} \big( q - q^{-1} \big)^j \, {[j]}_q!
\left\{ {d \, ; s} \atop {m \, , j} \right\} \, \c^{(j)} \, \b^{(j)}
\left\{ {a \, ; r} \atop {n \, , j} \right\}  \cr
   {{a \, ; r} \choose t} \, \b^{(n)} \, = \, \b^{(n)} {{a \, ; r+n}
\choose t} \quad ,  \qquad  {{a \, ; r} \choose t} \, \c^{(n)} \, =
\, \c^{(n)} {{a \, ; r+n} \choose t}  \cr
   {{d \, ; s} \choose t} \, \b^{(n)} \, = \, \b^{(n)} {{d \, ; s-n}
\choose t} \quad ,  \qquad  {{d \, ; s} \choose t} \, \c^{(n)} \, =
\, \c^{(n)} {{d \, ; s-n} \choose t} \; .  \cr }  $$
   \indent   Moreover, $ \calfqm $  is a  $ \Zqqm $--bialgebra,
whose bialgebra structure is given by
  $$  \displaylines{
   \Delta \left( \!\! {{a \, ; r} \choose n} \!\! \right) \, = \;
{\textstyle \sum\limits_{k=0}^n} \; q^{k (r-n)} \big( q - q^{-1} \big)^k
\, [k]_q! \cdot \Big( \b^{(k)} \otimes 1 \Big) \cdot \left\{\!
{a \otimes a \, ; \, r - k} \atop {n \, , \, k} \!\right\}
\cdot \Big( 1 \otimes \c^{(k)} \Big)  \cr
   \Delta \left( \b^{(n)} \right) \, = \, {\textstyle \sum\limits_{k=0}^n}
\; q^{-k(n-k)} \cdot a^k \, \b^{(n-k)} \otimes \b^{(k)} \, d^{n-k}  \cr
   \Delta \left( \c^{(n)} \right) \, = \, {\textstyle \sum\limits_{k=0}^n}
\; q^{-k(n-k)} \cdot \c^{(k)} \, d^{n-k} \, \otimes a^k \, \c^{(n-k)}  \cr
   \Delta \left( \!\! {{d \, ; s} \choose m} \!\! \right) \, = \;
{\textstyle \sum\limits_{k=0}^m} \; q^{k(s-m)} \big( q - q^{-1} \big)^k
\, [k]_q! \cdot \Big( 1 \otimes \b^{(k)} \Big) \cdot \left\{ {d \otimes d
\, ; \, s - k} \atop {m \, , \, k} \right\} \cdot \Big( \c^{(k)} \otimes 1
\Big)  \cr
   \hfill   \epsilon \left( \!\!
{{a \, ; r} \choose n} \!\! \right) = {r \choose n}_{\!q} \, ,  \hskip9pt
\epsilon \left( \b^{(\ell)} \right) = 0 \, ,  \hskip9pt
\epsilon \left( \c^{(\ell)} \right) = 0 \, ,  \hskip9pt
\epsilon \left( \!\! {{d \, ; s} \choose m} \!\! \right) =
{s \choose m}_{\!q}   \hfill   (\, \ell, m, n \geq 1)  \cr }  $$
(notation of \S\; 1.2) where  $ \; a = 1 + (q-1) \Big(\! {{a \, ; \, 0}
\atop 1} \!\Big) \, $,  $ \; d = 1 + (q-1) \Big(\! {{d \, ; \, 0}
\atop 1} \!\Big) \, $,  \; and terms like  $ \; \Big(\! {{x \otimes x
\, ; \, \sigma} \atop t} \!\Big) \; $  (with  $ \, x \in \{a,d\,\} \, $,
$ \, \sigma \in \Z \, $  and  $ \, t \in \N_+ \, $)  must be expanded
following the rule
  $$  \displaylines{
   {{x \otimes x \, ; 2 \tau} \choose t} \, = \, \sum_{r+s=\nu} q^{-sr}
{{x \, ; \tau} \choose r} \! \otimes {{x \, ; \tau} \choose s} \, x^r
\, = \, \sum_{r+s=\nu} q^{-rs} x^s {{x \, ; \tau} \choose r} \! \otimes
{{x \, ; \tau} \choose s}  \cr
   {{x \otimes x \, ; 2 \tau \! + \! 1} \choose t} \! = \!\!\!
\sum_{r+s=\nu} \!\!\! q^{-(1-s)r} {{x \, ; \tau} \choose r} \otimes
{{x \, ; \tau \! + \! 1} \choose s} \, x^r = \!\!\! \sum_{r+s=\nu}
\!\!\! q^{-(1-r)s} x^s {{x \, ; \tau \! + \! 1} \choose r} \!
\otimes {{x \, ; \tau} \choose s}  \cr }  $$
according to whether  $ \sigma $  is even  ($ \, = 2 \tau $)  or
odd  ($ \, = 2 \tau + 1 $),  and consequently for  $ \; \left\{\!
{x \otimes x \, ; \, \sigma} \atop {t \, , \, \ell} \!\!\right\} \, $.
                                     \hfill\break
   \indent   In particular,  $ \calfqm $  is a  $ \Zqqm $--integral
form (as a bialgebra) of  $ \fqm \, $.
 \vskip4pt
   (b) \;  $ \calfqgl $  is the unital associative  $ \Zqqm $--algebra
with generators
  $$  {{a \, ; r} \choose n} \; ,  \qquad  \b^{(n)} \; ,  \qquad
\c^{(n)} \; ,  \qquad  {{d \, ; r} \choose n} \; ,  \qquad
{D_q}^{\!-1}   \eqno \forall \;\;\, n \in \N \, , \, r \in \Z
\hskip25pt {}  $$
and relations as in claim (a), plus the additional relations
  $$  \displaylines{
   \hfill   {D_q}^{\!-1} \, {{x \, ; r} \choose n} \, = \,
{{x \, ; r} \choose n} \, {D_q}^{\!-1} \; ,  \quad  {D_q}^{\!-1}
\, \hbox{\bf y}^{(n)} \, = \, \hbox{\bf y}^{(n)} \, {D_q}^{\!-1}
\hfill \big(\, \hbox{\bf y} \in \big\{ \b, \c \big\} ,
\, r \in \Z \, ,  \, n \in \N \,\big)  \cr
   {D_q}^{\!-1} \, + \, (q-1) \, {{a \, ; 0} \choose 1}
\, {D_q}^{\!-1} \, + \, (q-1) {{d \, ; 0} \choose 1} \,
{D_q}^{\!-1} \, +   \hfill  \cr
   \hfill   + \; {(q-1)}^2 {{a \, ; 0} \choose 1} {{d \, ; 0}
\choose 1} \, {D_q}^{\!-1} \, - \, q \, \big( q - q^{-1} \big)^2
\, \b^{(1)} \c^{(1)} \, {D_q}^{\!-1} \; = \; 1 \;\; .  \cr }  $$
   \indent   Moreover, $ \calfqgl $  is a Hopf  $ \, \Zqqm $--algebra,
whose Hopf algebra structure is given by the same formul\ae{}  as in
claim (a) for  $ \Delta $  and  $ \epsilon $  plus the formul\ae{}
  $$  \displaylines{
   S \left( \!\! {{a \, ; r} \choose n} \!\! \right) = \, {{d \,
{D_q}^{\!-1} \, ; \, r} \choose n} \; ,  \qquad \qquad  S \left(
\b^{(n)} \right) = \, {(-1})^n q^{-n} \, \b^{(n)} \, {D_q}^{\!-n}  \cr
   \phantom{\; .}  S \left( \c^{(n)} \right) = \, {(-1})^n q^{+n}
\c^{(n)} \, {D_q}^{\!-n} \; ,  \qquad \qquad  S \left( \!\!
{{d \, ; r} \choose n} \!\! \right) = \, {{a \, {D_q}^{\!-1}
\, ; \, r} \choose n}
%
%
%
     \cr
   \Delta\big({D_q}^{\!-1}\big) \, = \, {D_q}^{\!-1} \otimes
{D_q}^{\!-1} \; ,  \qquad  \varepsilon\big({D_q}^{\!-1}\big) \, =
\, 1 \; ,  \qquad  S\big({D_q}^{\!-1}\big) \, = \, D_q \; .  \cr }  $$
   \indent   In particular,  $ \calfqgl $  is a  $ \Zqqm $--integral
form (as a Hopf algebra) of  $ \fqgl \, $.
 \vskip4pt
   (c) \;  $ \calfqsl $  is generated, as a unital associative
$ \Zqqm $--algebra,  by generators as in claim (a).  These generators
enjoy all relations in (a), plus some additional relations, springing
out of the relation  $ \, D_q = 1 \, $  in  $ \fqsl \, $.  Moreover,
$ \calfqsl $  is a Hopf  $ \, \Zqqm $--algebra,  whose Hopf algebra
structure is given as in (a) and (b), but setting  $ \, D_q = 1 \, $.
                                                   \par
   In particular,  $ \calfqsl $  is a  $ \, \Zqqm $--integral
form (as a Hopf algebra) of  $ \fqsl \, $.
\endproclaim

\demo{Proof}  {\it (a)} \,  Thanks to  Theorem 3.1{\it (a)},  the set of
elements considered in the statement does generate  $ \calfqm \, $.
%
%
 As for relations, the third line ones are those
springing out of the relation  $ \; a \, d - d \, a \, = \, \big(
q - q^{-1} \big) \, b \, c \; $  in  $ \fqm $;  those in fourth and
fifth line are the ones following from the relations  $ \; a \, b
\, = \, q \, b \, a \; $,  $ \; a \, c \, = q \, d \, c \; $,  $ \;
b \, d \, = \, q \, d \, b \; $  and  $ \; c \, d \, = q \, d \, c
\; $;  \; the first line ones follow from  $ \; b \, c \, = \, c
\, b \; $,  \; and those in second line are obvious.  The sole
non-trivial relations are the third line ones, which we shall
now prove, by induction on  $ m \, $.
                                          \par
   We set  $ \, A^r_{n,k} := \left\{\! {a \, ; \,\, r} \atop {n \, , \, k}
\!\right\} \, $,  $ \, D^s_{m,h} := \left\{\! {d \, ; \, s} \atop
{m \, , \, h} \!\right\} \, $.  The basis of induction  ($ m = 1 $)
follows from
  $$  \big[ A^t_{n,k} \, , d \,\big] \, = \; q^{t+k-n} \big( q -
q^{-1} \big)^2 \, \c^{(1)} \, \b^{(1)} A^t_{n,k+1}   \eqno (3.4)  $$
which in turn is easily proved by induction on  $ k $  using the
commutation formula
  $$  \textstyle{ \Big[ \Big(\! {{a \, ; \, \ell \,} \atop u} \!\Big) \, ,
\, d \, \Big] } \, = \, q^{\ell - u} \big( q - q^{-1} \big)^2 \,
\c^{(1)} \, \b^{(1)} A^\ell_{u,1}   \eqno (3.5)  $$
%
%
that directly follows from definitions and the relation  $ \; a \, d
- d \, a \, = \, \big( q - q^{-1} \big) \, b \, c \; $.
                                        \par
   For the general case  ($ \, m > 1 $),  using formul\ae{}  (3.4--5)
we get
  $$  \displaylines{
  \left[ \textstyle{ \Big(\! {{a \, ; \, r} \atop n} \!\Big) } ,
\textstyle{ \Big(\! {{d \, ; \, s} \atop {m+1}} \!\Big) } \right]
= \left[ \textstyle{ \Big(\! {{a \, ; \, r} \atop n} \!\Big) },
\textstyle{ \Big(\! {{d \, ; \, s} \atop m} \!\Big) } \right] \,
\textstyle{ {{\; d \, q^{s-m} - 1 \;} \over {\; q^{m+1} - 1 \;}} }
\, + \, \textstyle{ \Big(\! {{d \, ; \, s} \atop m} \!\Big) }
\left[ \textstyle{ \Big(\! {{a \, ; \, r} \atop n} \!\Big) } \, ,
\textstyle{ {{\; d \, q^{s-m} - 1 \;} \over {\, q^{m+1} - 1 \,}} } \right] =   \hfill  \cr
%
%
   = \big( q^{m+1} - 1 \big)^{-1} \bigg( {\textstyle
\sum\limits_{j=1}^{n \wedge m}} \, q^{j((r+s)-(n+m))} \big( q - q^{-1}
\big)^j q^{{j \choose 2}} \, {[j]}_q! \, D^s_{m,j} \left( q^{s-k+2j} d
- 1 \right) \c^{(j)} \, \b^{(j)} A^r_{n,j} \, +  \cr
   \hfill   + \, {\textstyle \sum_{i=1}^{n \wedge m}} \,
q^{i \, ((r+s)-(n+m)) + {i \choose 2}} \big( q \! - \! q^{-1} \big)^{i+1}
\, {[i\!+\!1]}_q! \, q^{r+i+s} \big( q^i \! - \! q^{-i} \big) \, D^s_{m,i}
\, \c^{(i+1)} \, \b^{(i+1)} A^r_{n,i+1} \bigg) }  $$
   Comparing the previous result with the expected formula for
$ (m+1) $,  we see that the latter holds if and only if the
following identity holds
  $$  q^j \, D^s_{m,j} \, \big( q^{s-m+2j} d - 1 \big) \, + \,
\big( q^{2j} - 1 \big) \, D^s_{m,j-1} \; = \; \big( q^{m+1} - 1 \big)
D^s_{m+1,j}  $$
which is just a special case of  Lemma 4.5{\it (a)\/}  later on.
                                        \par
   The previous analysis shows that the given relations do hold in
$ \calfqm \, $;  \, in order to prove the claim, we must show that
these generate the ideal of all possible relations.  This amounts
to show that the algebra enjoying only the given relations is in
fact isomorphic to  $ \calfqm \, $.  To this end, it is enough to
prove the following.  Let  $ \Cal{B}' $  be any one of the PBW--like
bases provided by Theorem 3.1{\it (a)\/}:  then the given relations
are enough to expand any product of the given generators as a
$ \Zqqm $--linear  combination of the monomials in  $ \Cal{B}' \, $.
                                        \par
   Now, if  $ \, \Cal{B}' = \left\{\, \Big(\! {{d \, ; \, 0} \atop
\delta} \!\Big) \, \c^{(\kappa)} \; \b^{(\beta)} \, \Big(\! {{a \, ; \, 0}
\atop \alpha} \!\Big) \, \;\Big|\; \alpha, \beta, \kappa, \delta \in \N
\,\right\} \, $,  \,
%
%
 then the given relations clearly allow to write any product of the
generators as an element of the  $ \Zqqm $--span  of  $ \Cal{B}' \, $.
                                        \par
   As to the bialgebra structure, everything is just a matter of
computation.  Yet we can point out just one key detail: namely, by
definition we have
 \vskip-9pt
  $$  \textstyle{ \Delta \left( \! \Big(\! {{a \, ; \, r} \atop n} \!\Big)
\! \right) \, = \; \Big( {{\Delta(a) \, ; \, r} \atop n} \Big) \, = \;
\Big( {{a \, \otimes \, a \, + \, b \, \otimes \, c \; ; \;\; r} \atop n} \Big)
\, = \; \bigg( {{a \, \otimes \, a \, + \, \big( q - q^{-1} \big)^2 \;
\text{\bf b} \, \otimes \, \text{\bf c} \,\; ; \;\; r} \atop
\phantom{\big|} n \phantom{\big|}} \bigg) }  $$
 \vskip-3pt
\noindent
 and then one gets the formula in the claim via Lemma 4.2; similarly
for  $ \, \Delta \left(\! \Big(\! {{d \, ; \, s} \atop m} \!\Big)
\!\right) \, $.
 \vskip4pt
   {\it (b)} \,  The fact that  $ \calfqgl $  admits the given presentation
is a direct consequence of claim  {\it (a)\/}  and of  Theorem 3.1{\it (b)},
but for the additional relations in first line.  The latter mean that
$ {D_q}^{-1} $  is central (because  $ D_q $  is) while the second line relation
is a reformulation of the relation  $ \, D_q \, {D_q}^{\!-1} = 1 \, $.  The statement on the Hopf structure also follows from claim  {\it (a)\/}  and
Theorem 3.1{\it (b)\/}  and from the  formul{\ae}  for the antipode in
$ \fqgl $  (cf.~\S\; 2.2), but for the  formul{\ae}  for  $ {D_q}^{\!-1} $
which follow from  $ \; \Delta(D_q) = D_q \otimes D_q \, $,  $ \;
\varepsilon(D_q) = 1 \, $,  $ \; S\big(D_q\big) = {D_q}^{\!-1} \; $.
%
%
                                            \par
   Now, the  formul{\ae}  for  $ \Delta $  and  $ \epsilon $
show that  $ \calfqgl $  is a  $ \Zqqm $--subbialgebra
of  $ \fqgl \, $.  For the antipode,  $ \, \Big\langle S \big(
\calfqgl \big) \, , \, \caluqgl \Big\rangle = \Big\langle \calfqgl
\, , \, S\big(\caluqgl\big) \! \Big\rangle \subseteq \Zqqm \, $,  \,
which gives  $ \, S \big( \calfqgl \big) \subseteq \calfqgl \; $.
The claim follows.
 \vskip4pt
   {\it (c)} \,  This follows again from claim  {\it (a)\/}  and
Theorem 3.1.   \qed
\enddemo

\vskip9pt

\proclaim{Corollary 3.3}  For every  $ \, X \in \{M, GL\} \, $,  \,
let  $ \big( D_q - 1 \big) $  be the two-sided ideal of  $ F_q[X_2] $
generated by  $ \big( D_q - 1 \big) $,  and let  $ \, \Cal{D}(X_2) :=
\big( D_q - 1 \big) \cap \calF_q[X_2] \, $,  \, a two-sided ideal
of  $ \calF_q[X_2] \, $.
                                                \par
   \indent (a) \; The epimorphism  $ \; \fqm \,{\buildrel \pi \over
{\llongtwoheadrightarrow}}\, \fqm \Big/ \big( D_q - 1 \big) \, \cong
\, \fqsl \; $  restricts to an epimorphism  $ \; \calfqm \,{\buildrel
\pi \over {\llongtwoheadrightarrow}}\, \calfqm \Big/ \Cal{D}(M_2) \,
\cong \, \calfqsl \; $  of  $ \, \Zqqm $--bialgebras.
                                                \par
   \indent (b) \; The epimorphism  $ \; \fqgl \,{\buildrel \pi \over
{\llongtwoheadrightarrow}}\, \fqgl \Big/ \big( D_q - 1 \big) \, \cong
\, \fqsl \; $  restricts to an epimorphism  $ \; \calfqgl \,{\buildrel
\pi \over {\llongtwoheadrightarrow}}\, \calfqgl \Big/ \Cal{D}(GL_2)
\, \cong \, \calfqsl \; $  of Hopf  $ \, \Zqqm $--algebras.   \qed
\endproclaim

\vskip7pt

  {\bf 3.4 Remarks:} \,
%
%
 {\it (a)} \, Besides those given
in Theorem 3.2, there are several other (equivalent) commutation relations
between the  $ \Big(\! {{a \, ; \, r} \atop n} \!\Big) $'s  and the
$ \Big(\! {{d \, ; \, s} \atop m} \!\Big) $'s   ---   see [GR1],
\S\; 3.3, for details.
                                                   \par
    {\it (b)} \, One should compute an expression for  $ \, S \Big(\!
\Big(\! {{a \, ; \, r} \atop n} \!\Big) \!\Big) = \Big(\! {{d \,
{D_q}^{\!-1} \, ; \; r} \atop n} \!\Big) \, $  and  $ \, S \Big(\! \Big(\!
{{d \, ; \, r} \atop n} \!\Big) \!\Big) = \Big(\! {{a \, {D_q}^{\!-1} \, ;
\; r} \atop n} \!\Big) \, $  in terms of the generators in Theorem 3.2!
%
%
 In fact, this is a very tough task.
%
%

\vskip11pt

   {\bf 3.5 Relations in  $ \calfqsl \, $.} \, By  Theorem 3.2{\it (c)},  in
$ \calfqsl $  the generators  $ \Big( {{a \, ; \, r} \atop n} \Big) \, $,
$ \, \b^{(n)} \, $,  $ \, \c^{(n)} \, $  and  $ \, \Big( {{d \, ; \, r} \atop n} \Big) \, $   --- for all  $ \, n \in \N \, $,  $ \, r \in \Z \, $  ---   enjoy
the same relations as in Theorem 3.2  {\sl plus some additional ones},  springing
out of the relation  $ \; a \, d - q \, b\, c \, = \, 1 \; $  in  $ \fqsl \, $.
 \vskip3pt
   A  {\bf first (set of) relation(s)}  is the following (for
all  $ \, \alpha, \delta \in \N \, $,  $ \, r, s \in \Z \, $):
  $$  \displaylines{
   {(\alpha)}_q \, {(\delta)}_q \cdot \textstyle{ \Big(\! {{a \, ; \, r}
\atop \alpha} \!\Big)} \, \textstyle{ \Big(\! {{d \, ; \, s} \atop \delta}
\!\Big) } \, + \, {(\alpha)}_q \cdot \textstyle{ \Big(\! {{a \, ; \, r}
\atop \alpha} \!\Big) } \, \textstyle{ \Big(\! {{d \, ; \, s} \atop
{\delta - 1}} \!\Big) } \, +
\hfill  \cr
    \hfill   + \, {(\delta)}_q \cdot \textstyle{ \Big(\! {{a \, ; \, r} \atop
{\alpha - 1}} \!\Big) } \, \textstyle{ \Big(\! {{d \, ; \, s} \atop \delta}
\!\Big) } \, + \, {(2 - \alpha - \delta + r + s)}_q \cdot \textstyle{ \Big(\!
{{a \, ; \, r} \atop {\alpha - 1}} \!\Big) } \, \textstyle{ \Big(\!
{{d \, ; \, s} \atop {\delta - 1}} \!\Big) } \, +   \hfill
%
%
 \cr
   \hfill   + \; q^{1 - \alpha - \delta + r + s} \, (q \! - \! 1) \,
{{(2)}_q}^{\!\! 2} \cdot \textstyle{ \Big(\! {{a \, ; \, r} \atop
{\alpha - 1}} \!\Big) } \, \b^{(1)} \, \c^{(1)} \,
\textstyle{ \Big(\! {{d \, ; \, s} \atop {\delta - 1}} \!\Big) }
\; = \; 0  \quad \qquad  \cr }  $$
%
%

\vskip4pt

   A  {\bf second (set of) relation(s)}  is (for
all  $ \, h, n \in \N \, $,  $ \, r, s \in \Z \, $):
  $$  \hskip-7pt   \hbox{ $ \eqalign{
   \underset{(j,\ell) \not= (0,0)}\to{{\textstyle \sum\limits_{j=0}^h}
\, {\textstyle \sum\limits_{\ell=0}^n}}
\! q^{{j \choose 2} + {\ell \choose 2}} \, {(q-1)}^{j+\ell-1} \, {(j)}_q!
\, {(\ell)}_q! \, \textstyle{ \Big(\! {h \atop j} \!\Big) }_{\!q} \,
\textstyle{ \Big(\! {n \atop \ell} \!\Big) }_{\!q} \cdot
\textstyle{ \Big(\! {{a \, ; \, r} \atop j} \!\Big) } \,
\textstyle{ \Big(\! {{d \, ; \, s} \atop \ell} \!\Big) }
\,  &  =  \; {\big( n(r+s) \big)}_q \, +  \cr
   + \; q^{n(r+s)} \! \underset{(j,\ell) \not= (0,0)}
\to{{\textstyle \sum\limits_{i=0}^{h-n}} {\textstyle \sum\limits_{t=0}^n}}
\! q^{t^2 + {i \choose 2}} \, {(q\!-\!1)}^{i-1} \, {\big( q \! - \! q^{-1}
\big)}^{\!2t} \, {(i)}_q! \, {{[t]}_q!}^2 \,
\textstyle{ \Big(\! {{h-n} \atop i} \!\Big) }_{\!\!q}  &
\textstyle{ \Big(\! {n \atop t} \!\Big) }_{\!\!q^2} \!
\cdot \textstyle{ \Big(\! {{a \, ; \, r} \atop i} \!\Big) }
\, \b^{(t)} \, \c^{(t)}  \cr } $ } %
%
 $$
  $$  \hskip-7pt  \hbox{ $ \eqalign{
   \underset{(j,\ell) \not= (0,0)}\to{{\textstyle \sum\limits_{j=0}^h}
\, {\textstyle \sum\limits_{\ell=0}^n}}
\! q^{{j \choose 2} + {\ell \choose 2}} \, {(q-1)}^{j+\ell-1} \, {(j)}_q!
\, {(\ell)}_q! \, \textstyle{ \Big(\! {h \atop j} \!\Big) }_{\!q} \,
\textstyle{ \Big(\! {n \atop \ell} \!\Big) }_{\!q}
\cdot \textstyle{ \Big(\! {{a \, ; \, r} \atop j} \!\Big) } \,
\textstyle{ \Big(\! {{d \, ; \, s} \atop \ell} \!\Big) } \,  &  =
\; {\big(h(r+s)\big)}_q \, +  \cr
   + \; q^{h(r+s)} \! \underset{(i,t) \not= (0,0)}
\to{{\textstyle \sum\limits_{i=0}^{n-h}} {\textstyle \sum\limits_{t=0}^h}}
\! q^{t^2 + {i \choose 2}} \, {(q\!-\!1)}^{i-1} \, {\big( q \! - \! q^{-1}
\big)}^{\!2t} \, {(i)}_q! \, {{[t]}_q!}^2 \,
\textstyle{ \Big(\! {{n-h} \atop i} \!\Big) }_{\!q}  &
\textstyle{ \Big(\! {h \atop t} \!\Big) }_{\!q^2} \!
\cdot \, \b^{(t)} \, \c^{(t)} \,
\textstyle{ \Big(\! {{d \, ; \, s} \atop i} \!\Big) }  \cr } $ }
%
%
 $$
where
%
%
 the first identity
holds for all  $ \, h \geq n \, $
          and
%
%
%
 the second
for all  $ \, h \leq n \, $.

\vskip4pt

   A  {\bf third (set of)
       relation(s)}  is
%
 for all  $ \, n \in \N_+ \, $,
%
%
  $$  \textstyle{ \Big(\! {{a \, ; \, 0} \atop n} \!\Big) } +
\textstyle{ \Big(\! {{d \, ; \, 0} \atop n} \!\Big) }  \; = \;
{\textstyle \sum_{h=1}^n} \, {\widetilde\alpha}_h^{\,n} \,
\b^{(h)} \, \c^{(h)} \, + \, {\textstyle \sum_{\Sb  h=1  \\
h \leq i  \endSb}^n} \, {\widetilde\beta}_{h,i}^{\,n} \,
\left(\! \textstyle{ \Big(\! {{a \, ; \, 0} \atop h} \!\Big) }
\textstyle{ \Big(\! {{d \, ; \, 0} \atop i} \!\Big) } +
\textstyle{ \Big(\! {{a \, ; \, 0} \atop i} \!\Big) }
\textstyle{ \Big(\! {{d \, ; \, 0} \atop h} \!\Big) } \!\right)
%
%
 $$
with (for all  $ \, 1 \leq h \leq i \leq n \, $)
  $$  \eqalign{
   {\widetilde\alpha}_h^{\,n} \,  &  = \; q^{h^2 + {{n-h+1} \choose 2}}
\, {{[h]}_q!}^2 \,
\textstyle{
  {{\, {(q - q^{-1})}^{2h}\,} \over {\, {(n)}_q! \, {(q - 1)}^h \,}}
\; \, {\Big(\! {n \atop {2(n-h)}} \!\Big)}_{\!q} }
{\big(2(n\!-\!h)\!-\!1\big)}_q!!  \cr
   {\widetilde\beta}_{h,i}^{\,n} \,  &  = \; -
{(1 + \delta_{h,i})}^{-1} \, q^{{n-h-i+1} \choose 2} \, {(q-1)}^{h+i-n} \,
\textstyle{
 {{\,{(h)}_q!\,} \over {\,{(n-i)}_q!\,}} \,
{\Big(\! {i \atop {n-h}} \!\Big)}_{\!q} }
 \cr }  $$
   \indent   All these  formul{\ae}  are proved in detail in [GR1].

\vskip10pt

   Further byproducts of Theorem 3.2 concern the specializations of
$ \calfqm \, $,  $ \calfqgl $  and  $ \calfqsl $  at roots of unity,
including the case  $ \, q = 1 \, $,  \, as follows:

\vskip11pt

\proclaim{Corollary 3.6}
 \vskip3pt
   {\it (a)} \,  There exists a  $ \Z $--bialgebra isomorphism
$ \; \calF_1[M_2] \cong U_\Z\big({\gergl_2}^{\!*} \big) \, $  given by
 \vskip-13pt
%
%
  $$  \textstyle { \Big(\! {{a \, ; \, 0} \atop n} \!\Big){\Big|}_{q=1} \!\!
\mapsto \Big(\! {\g_1 \atop n} \!\Big) \; ,  \quad  \b^{(n)}{\Big|}_{q=1} \!\!
\mapsto \f^{\,(n)} \; ,  \quad  \c^{(n)}{\Big|}_{q=1} \!\! \mapsto \e^{(n)} \; ,
\quad  \Big(\! {{d \, ; \, 0} \atop m} \!\Big){\Big|}_{q=1} \!\! \mapsto
\Big(\! {\g_2 \atop m} \!\Big) \quad . }  $$
 \vskip-3pt
\noindent
 In particular  $ \, \calF_1[M_2] \, $  is a Hopf  $ \, \Z $--algebra,
isomorphic to  $ \, U_\Z\big({\gergl_2}^{\!\!*} \big) \, $.
 \vskip3pt
   {\it (b)} \,  There exists a Hopf  $ \, \Z $--algebra
isomorphism  $ \, \calF_1[{GL}_2] \cong U_\Z\big({\gergl_2}^{\!\!*}
\big) \, $,  \, which is uniquely determined by the formul{\ae}
in claim (a).
 \vskip3pt
   {\it (c)} \,  There exists a Hopf  $ \, \Z $--algebra  isomorphism
$ \,\calF_1[{SL}_2] \cong U_\Z\big({\gersl_2}^{\!\!*} \big) $  given
by the same formul{\ae}  as in claim (a), where one must read  $ \,
\g_1 = \h \, $,  $ \, \g_2 = -\h \, $.   \qed
\endproclaim

\demo{Proof} At  $ \, q = 1 \, $,  \, Theorem 3.2{\it(a)\/}  provides a presentation for  $ \calF_1[M_2] \, $;  a straightforward comparison then shows that the latter
is the standard presentation of  $ \, U_\Z\big({\gergl_2}^{\!\!*} \big) \, $,
following the correspondence given in the claim (actually, in the first presentation one also has the specialization at  $ \, q = 1 \, $  of the  $ \Big(\! {{x \, ; \, r} \atop k} \!\Big) $'s   --- with  $ \, x \in \{a,d\} \, $  ---   but these are
generated by the specializations of the  $ \Big(\! {{x \, ; \, 0} \atop \nu} \!\Big) $'s).  This yields a  $ \Z $--algebra  isomorphism: a moment's check shows that it is one of  $ \Z $--bialgebras  too.  This proves {\it (a)\/}  and,    
   \hbox{with minimal changes,  {\it (b)\/}  and  {\it (c)\/}  too.   \qed}
\enddemo

\vskip11pt

\proclaim{Proposition 3.7}  Let  $ \varepsilon $  be a root of unity,
of odd order, and apply notation of\/ \S\, 1.3.
 \vskip3pt
   (a) \,  The specialization  $ \, \calfem \longhookrightarrow
\gerH_\varepsilon^g \; $  at  $ \, q = \varepsilon \, $  of the embedding
$ \, \calfqm \longhookrightarrow \gerH_q^g \; $  is a  $ \, \Zeps $--algebra  isomorphism.
%
%
 \vskip3pt
   (b) \,  The embedding  $ \, \calfem \longhookrightarrow \calfegl \, $
of  $ \, \Zeps $--bial\-gebras  is an isomorphism.  In particular,
$ \calfem $  and  $ \gerH^g_\varepsilon $  both are Hopf
$ \Zeps $--algebras  isomorphic to  $ \calfegl \, $.
 \vskip3pt
   (c) \,  The specialization  $ \calfesl $  is a Hopf  $ \, \Zeps $--algebra,  isomorphic to  $ \gerH_\varepsilon^s $  via the specialization of the embedding
$ \, \calfqsl \longhookrightarrow \gerH_q^s \;\, $.
\endproclaim

\demo{Proof} The embedding  $ \; \calfqm \longhookrightarrow \! \calfqgl \; $
induces an embedding  $ \; \calfem \longhookrightarrow \! \calfegl \, $.
Then  {\it (a)\/}  will follow by proving that  $ \, D_\varepsilon := D_q
\mod (q - \varepsilon) \, \calfqm \, $  is invertible in  $ \calfem \, $.
                                                 \par
   Let  $ \ell $  be the (multiplicative) order of  $ \varepsilon \, $.
Lemma 4.3 gives  $ \, a^\ell = 1 \, $  in  $ \calfem \, $,  \, so  $ a $  is invertible in  $ \calfem $  with  $ \, a^{-1} = a^{\ell-1} \in \calfem \, $,
\, and also  $ \, \b^{\,\ell} = 0 \in \calfem \, $,  \, so  $ \, b^{\,\ell}
= {\big( \varepsilon - \varepsilon^{-1} \big)}^\ell \, \b^{\,\ell} = 0 \, $
in  $ \calfem \, $.  Similarly,  $ \, d^{\,-1} = d^{\,\ell-1} \in \calfem
\, $  and  $ \, c^{\,\ell} = \, 0 \, $.  The power series expansion of  $ \,
{(1-x)}^{-1} \, $  then gives
                                                       \par
   \centerline{ $ \; {D_\varepsilon}^{\!-1} \, = \,
%
%
 \big(\, 1 - \varepsilon \, b \, d^{-1}
a^{-1} \, c \,\big)^{-1} d^{-1} a^{-1} \,
= \, \sum\limits_{n=0}^{\ell-1} \varepsilon^n
\, b^n \big( d^{\,\ell-1} a^{\ell-1} \big)^n c^n \,
d^{\,\ell-1} a^{\ell-1} \in \calfem \;\; $. }
%
%
%
%
%
                                               \par
   As to the second part, note that the embedding  $ \, \xi : \fqgl
\longhookrightarrow \H^g_q \, $  extends to an identity  $ \, \fqgl
\big[ d^{-1} \big] \! = \H^g_q \, $:  \, this comes from [DL], \S\;
1.8 (adapted to the case of  $ {GL}_2 $)  or directly from the explicit
description of  $ \xi $  in \S\; 2.4.  This yields also  $ \, \gerH^g_q
= \calfqgl \big[ d^{-1} \big] \, $.  In fact, if  $ \, \eta \in \gerH^g_q
= \calfqgl \big[ d^{-1} \big] \, $  then there is  $ \, n \in \N \, $  such
that  $ \, \eta \, d^n \in \calfqgl \, $,  \, and also  $ \; \big\langle \eta
\, d^n , \, \caluqgl \big\rangle = \Big\langle \eta \otimes d^{\otimes n} ,
\, \Delta^{(n+1)} \big( \caluqgl \big) \Big\rangle \subseteq \big\langle
\eta , \, \caluqgl \big\rangle \cdot \big\langle d , \, \caluqgl
\big\rangle^n \subseteq \Zqqm \; $  because  $ \, \eta , d \in \gerH^g_q
\, $.  Thus  $ \; \eta \, d^n \in \fqgl \bigcap \gerH^g_q = \calfqgl
\, $,  \; whence  $ \; \eta \in
%
%
 \calfqgl \big[
d^{-1} \big] \, $;  \; the outcome is  $ \, \gerH^g_q \subseteq
\calfqgl \big[ d^{-1} \big] $,  \, and the converse is clear.  Now
$ \, \gerH^g_\varepsilon = \calfegl\big[d^{-1}\big] \, $;  \,
%
%
but we found  $ \, d^{-1} \in \calfem = \calfegl \, $,  \, so
$ \, \gerH^g_\varepsilon = \calfegl \big[ d^{-1} \big] \! =
\calfegl = \calfem \, $.
                                               \par
   The above proves claim  {\it (a)\/}  and  {\it (b)},  noting that
$ \, \calfegl \, $  is clearly a Hopf  $ \Zeps $--algebra.  Similarly,
{\it (c)\/}  can be proved like  {\it (a)},  or deduced from the
latter.   \qed
\enddemo

\vskip11pt

\proclaim{Theorem 3.8}  Let  $ \varepsilon $  be a root of unity, of
odd order  $ \ell \, $.
 \vskip3pt
   (a) \,  The quantum Frobenius morphism (2.1)
%
%
is defined over  $ \Zeps \, $,  i.e.{} it restricts to an epimorphism
of  $ \, \Zeps $--bialgebras
 $ \; {{\Cal F}r}_{M_2}^{\,\Z} \, \colon \, \calfem
\llongtwoheadrightarrow \, \Zeps \otimes_\Z \calF_1[M_2] \,
\cong \, \Zeps \otimes_\Z U_\Z\big({\gergl_2}^{\!*}\big) \; $,
\; coinciding, via Corollary 3.6 and Proposition 3.7, with (2.5),
and given on generators by
  $$  {\Cal{F}r}_{M_2}^{\,\Z} :
 \left\{ \hskip-5pt
    \hbox{ $ \matrix
   \Big(\! {a \, ; \, 0 \atop n} \!\Big) \!
\Big\vert_{q=\varepsilon} \hskip-7pt \mapsto
\Big( {a \, ; \, 0 \atop {n / \ell}} \Big) \!
\Big\vert_{q=1} \hskip-5pt = \Big( {\g_1 \atop {n / \ell}} \Big)
&  \hskip-6pt \hbox{ if } \; \ell \Big\vert n \; ,  &  \hskip-0pt
\Big(\! {a \, ; \, 0 \atop n} \!\Big) \! \Big\vert_{q=\varepsilon} \hskip-7pt
\mapsto \, 0  &  \hskip-6pt \hbox{ if } \; \ell \!\! \not\Big\vert n  \\
   \b^{(n)}\Big\vert_{q=\varepsilon} \hskip-1pt \mapsto \;
\b^{(n/\ell)}\Big\vert_{q=1} \hskip-3pt = \; \f^{\,(n/\ell)}
&  \hskip-6pt \hbox{ if } \; \ell \Big\vert n \; ,  &  \hskip-0pt
\b^{(n)}\Big\vert_{q=\varepsilon} \hskip-1pt \mapsto \; 0  &
\hskip-6pt \hbox{ if } \; \ell \!\! \not\Big\vert n  \\
   \hskip1,3pt \c^{(n)}\Big\vert_{q=\varepsilon} \hskip-1pt \mapsto
\;\hskip1,4pt
\c^{(n/\ell)}\Big\vert_{q=1} \hskip-3pt = \; \e^{(n/\ell)}
&  \hskip-6pt \hbox{ if } \; \ell \Big\vert n \; ,  &  \hskip-0pt
\c^{(n)}\Big\vert_{q=\varepsilon} \hskip-1pt \mapsto \; 0  &
\hskip-6pt \hbox{ if } \; \ell \!\! \not\Big\vert n  \\
   \Big(\! {d \, ; \, 0 \atop n} \!\Big) \!
\Big\vert_{q=\varepsilon} \hskip-7pt \mapsto
\Big( {d \, ; \, 0 \atop {n / \ell}} \Big) \!
\Big\vert_{q=1} \hskip-5pt = \Big( {\g_2 \atop {n / \ell}} \Big)
&  \hskip-6pt \hbox{ if } \; \ell \Big\vert n \; ,  &  \hskip-0pt
\Big(\! {d \, ; \, 0 \atop n} \!\Big) \! \Big\vert_{q=\varepsilon} \hskip-7pt
\mapsto \, 0  &  \hskip-6pt \hbox{ if } \; \ell \!\! \not\Big\vert n
            \endmatrix $ }
 \right.  $$
 \vskip3pt
   (b) \,  The quantum Frobenius morphism (2.2)
%
%
is defined over  $ \Zeps \, $,  i.e.{} it restricts to an epimorphism
of  $ \, \Zeps $--bialgebras
 $ \; {{\Cal F}r}_{{GL}_2}^{\,\Z} \, \colon \, \calfegl
\llongtwoheadrightarrow \, \Zeps \otimes_\Z \calF_1[{GL}_2] \,
\cong \, \Zeps \otimes_\Z U_\Z\big({\gergl_2}^{\!*}\big) \; $
\; coinciding, via Corollary 3.6 and Proposition 3.7, with (2.5) and
with  $ \, {{\Cal F}r}_{M_2}^{\,\Z} $  in (a).
 \vskip3pt
   (c) \,  The quantum Frobenius morphism (2.3)
%
%
is defined over  $ \Zeps \, $,  i.e.{} it restricts to an epimorphism
of  $ \, \Zeps $--bialgebras
 $ \; {{\Cal F}r}_{{SL}_2}^{\,\Z} \, \colon \, \calfesl
\llongtwoheadrightarrow \, \Zeps \otimes_\Z \calF_1[{SL}_2] \,
\cong \, \Zeps \otimes_\Z U_\Z\big({\gersl_2}^{\!*}\big) \; $
\; coinciding, via Corollary 3.6 and Proposition 3.7, with (2.4), and
described by  formul{\ae}  like in (a) with  $ \, \g_1 = +\h \, $  and
$ \, \g_2 = -\h \, $.
\endproclaim

\demo{Proof} By definition (cf.~[Ga1]) the morphism  $ \;
{{\Cal F}r}_{M_2}^{\,\Q} \, \colon \, \Qeps \otimes_{\Zeps} \calfem \longtwoheadrightarrow \, \Qeps \otimes_\Z U_\Z \big( {\gergl_2}^{\!*}
\big)
%
%
\; $  is the restriction (via  $ \, \widehat{\xi} : \calfqm \longhookrightarrow \gerH^g_q \, $  at  $ \, q = \varepsilon \, $  and  $ \, q = 1 \, $)  of the similar epimorphism  $ \; \gerFr_{{\gergl_2}^{\!*}}^{\,\Q} \, \colon \, \Qeps \otimes_{\Zeps} \gerH^g_\varepsilon \, \longtwoheadrightarrow \, \Qeps
\otimes_\Z \gerH^g_1
%
%
\; $  obtained by scalar extension from (2.5).  From this, direct computation
(taking into account that  $ \, {[\ell]}_\varepsilon
= 0 \, $)  gives, thanks to  Lemma 4.2{\it (a-1)},
  $$  \displaylines{
   {\frak F}{\frak r}_{{\gergl_2}^{\!*}}^{\,\Q}
\Big( {\widehat{\xi}\,{\Big|}_{q=\varepsilon}}
\Big( \! {\textstyle \Big(\! {{a \, ; \, 0} \atop n} \!\Big)}
\Big\vert_{q=\varepsilon} \Big) \Big) \, = \;
{\frak F}{\frak r}_{{\gergl_2}^{\!*}}^{\,\Z} \Big(
\! {\textstyle \Big(\! {{\Lambda_1 - \overline{F} \Lambda_2 \overline{E}
\; ; \; 0} \atop n} \!\Big)}\Big\vert_{q=\varepsilon} \Big) \, =
\hfill  \cr
   \hfill   = \, {\textstyle \sum\limits_{k=0}^{\ell - 1}} \;
\varepsilon^{- k n - {k \choose 2}} \big( \varepsilon^{-1} - \varepsilon
\big)^k {[k]}_\varepsilon! \, {\frak F}{\frak r}_{{\gergl_2}^{\!\!
*}}^{\,\Z} \! \left(\! F^{(k)} \right) \, {\frak F}{\frak
r}_{{\gergl_2}^{\!\!*}}^{\,\Z} \! \left(\! {\Lambda_2}^{\!k} \right) \,
{\frak F}{\frak r}_{{\gergl_2}^{\!\!*}}^{\,\Z} \! \left(\! E^{(k)} \right)
\, {\frak F}{\frak r}_{{\gergl_2}^{\!\!*}}^{\,\Z} \Big(\!\! \left(\!
{\textstyle \left\{ {{\Lambda_1 \, ; \; 0} \atop {n \; , \; k}} \right\}}
\right) \, =  \cr
   \hfill \hskip5pt   = \; {\frak F}{\frak r}_{{\gergl_2}^{\!\!*}}^{\,\Z}
\! \Big(\! {\textstyle \left\{ {{\Lambda_1 \, ; \; 0} \atop {n \; ,
\; 0}} \right\}}\Big\vert_{q=\varepsilon} \Big) \, = \;
{\frak F}{\frak r}_{{\gergl_2}^{\!\!*}}^{\,\Z} \! \Big(\!
{\textstyle \Big(\! {{\Lambda_1 \, ; \, 0} \atop n} \Big)}
\Big\vert_{q=\varepsilon} \Big) \, = \; \cases
     \hskip-3pt  \Big(\! {{\Lambda_1 \, ; \, 0} \atop {n / \ell}} \!\Big)
{\Big|}_{q=1}  \hskip-3pt  = \Big(\! {\g_1 \atop {n / \ell}} \!\Big)  &
\text{if} \quad \ell \,\big\vert\, n  \\
     \hskip-3pt  \quad\;  0  &  \text{if} \quad \ell \!\not\big\vert\, n
   \endcases  \cr }  $$
by the very description of  $ \gerFr_{{\gergl_2}^{\!*}}^{\,\Z} \, $;  \,
on the other hand, we also have
 $ \; {\widehat{\xi}\,{\Big|}_{q=1}} \Big( \! \Big(\! {{a \, ; \, 0}
\atop {n / \ell}} \Big)\Big\vert_{q=1} \Big) = \, \Big(\! {{\Lambda_1
- \overline{F} \Lambda_2 \overline{E} \; ; \; 0} \atop {n / \ell}} \Big)
\Big\vert_{q=1} \! = \, \left( \, \sum_{k=0}^{\ell - 1} \; q^{- k n -
{k \choose 2}} \, \big( q^{-1} \! - q \big)^k \, {[k]}_q! \cdot F^{(k)}
\, {\Lambda_2}^{\!k} \, E^{(k)} \, \left\{ {{\Lambda_1 \, ; \; 0} \atop
{n/\ell \; , \; k}} \right\} \right) \! \Big|_{q=1} \hskip-3pt =
 \break
\, \Big(\!
{{\Lambda_1 \, ; \, 0} \atop {n / \ell}} \Big)\Big\vert_{q=1} \hskip-1pt
= \, \Big(\! {\g_1 \atop {n / \ell}} \!\Big) \; $
by Corollary 3.6{\it (a)\/}  and  Lemma 4.2{\it (a-1)\/}  again.
This together with  $ \; \widehat{\xi}
{\,\big|}_{q=\varepsilon} \circ \big( \text{id}_{\Qeps}
\otimes_{\Zeps} \gerFr_{\gerg^*}^{\,\Z} \big) = \big(
\text{id}_{\Qeps} \otimes_{\Zeps} \widetilde{\,\xi}{\big|}_{q=1}
\big) \circ \calFr_G^{\,\Q} \; $  (see \S\; 2.4) give  $ \,
{\Cal{F}r}_{M_2}^{\,\Z} \Big( \! \Big(\! {{a \, ; \, 0} \atop n} \Big)
\Big\vert_{q=\varepsilon} \Big) = \Big(\! {{a \, ; \, 0} \atop
{n / \ell}} \Big) \Big\vert_{q=1} \hskip-3pt = \Big(\! {\g_1 \atop
{n / \ell}} \!\Big) \, $  if  $ \, \ell \big| n \, $,  \, and otherwise
$ \, {\Cal{F}r}_{M_2}^{\,\Z} \Big( \! \Big(\! {{a \, ; \, 0} \atop n}
\Big) \Big\vert_{q=\varepsilon} \Big) = 0 \, $,  \, as claimed.
Similarly,
 \vskip-5pt
  $$  \displaylines{
   {\frak F}{\frak r}_{{\gergl_2}^{\!*}}^{\,\Q} \Big( {\widehat{\xi}\,
\Big|_{q=\varepsilon}} \Big( \b^{(n)} \big|_{q=\varepsilon} \Big) \Big)
\, = \, {\frak F}{\frak r}_{{\gergl_2}^{\!*}}^{\,\Z} \Big( - F\big|_{q =
\varepsilon} \Lambda_2\big|_{q=\varepsilon} \Big)^{(n)} \, =
\hfill  \cr
   \hfill   = \, {(-1)}^n \varepsilon^{-{n \choose 2}} \cdot
{\frak F}{\frak r}_{{\gergl_2}^{\!*}}^{\,\Z} \big( F^{(n)} \big) \,
{\frak F}{\frak r}_{{\gergl_2}^{\!*}}^{\,\Z} \big( \Lambda^n \big) \,
= \, \cases
     \hskip-3pt  {(-1)}^n \varepsilon^{-{n \choose 2}} \cdot
F^{(n/\ell)}{\Big|}_{q=1} \hskip-3pt = \, {(-1)}^n
\varepsilon^{-{n \choose 2}} \cdot \text{f}^{\,(n/\ell)}  \\
     \hskip-3pt  \quad\;  0
     \endcases  \cr }  $$
 \vskip-2pt
\noindent
where the upper identity holds if  $ \, \ell \big| n \, $  and the
lower line holds if not.  On the other hand  $ \; {\widehat{\xi}\,
\Big|_{q=1}} \Big( \b^{(n/\ell)} \big|_{q =1} \Big) = {(-1)}^{n/\ell}
\varepsilon^{-{{n/\ell} \choose 2}} \cdot F^{(n/\ell)}{\Big|}_{q=1}
\hskip-3pt = \, {(-1)}^{n/\ell} \varepsilon^{-{{n/\ell} \choose 2}}
\cdot \text{f}^{\,(n/\ell)} \; $,  \; and a moment's check shows
that  $ \, {(-1)}^n \varepsilon^{-{n \choose 2}} = {(-1)}^{n/\ell}
\varepsilon^{-{{n/\ell} \choose 2}} \, $  whence we conclude.
Similarly
 \vskip-13pt
  $$  \displaylines{
   {\frak F}{\frak r}_{{\gergl_2}^{\!*}}^{\,\Q} \Big( {\widehat{\xi}\,
\Big|_{q=\varepsilon}} \Big( \c^{(n)} \big|_{q=\varepsilon} \Big) \Big)
\, = \, \cases
     \hskip-3pt  \varepsilon^{n \choose 2} \cdot
E^{(n/\ell)}{\Big|}_{q=1} \hskip-3pt = \,
\varepsilon^{n \choose 2} \cdot \text{e}^{(n/\ell)}  \\
     \hskip-3pt  \quad\;  0
     \endcases  \cr }  $$
 \vskip-3pt
\noindent
where the upper identity holds if  $ \, \ell \big| n \, $  and the
lower line if not, while  $ \; {\widehat{\xi}\,\Big|_{q=1}} \Big(
\c^{(n/\ell)} \big|_{q =1} \Big) = \, \varepsilon^{{n/\ell} \choose 2}
\cdot \e^{(n/\ell)} \, $,  \, and again one checks that  $ \,
\varepsilon^{n \choose 2} \! = \varepsilon^{{n/\ell} \choose 2}
\, $  so the claim for  $ \c^{(n)} $  follows.  Finally,
 \vskip-9pt
  $$  \textstyle
   {\frak F}{\frak r}_{{\gergl_2}^{\!*}}^{\,\Q}
\Big( {\widehat{\xi}\,{\Big|}_{q=\varepsilon}} \!
\Big( \! \Big(\! {{d \, ; \,0} \atop n} \Big) \Big\vert_{q=\varepsilon}
\Big) \Big) = \, {\frak F}{\frak r}_{{\gergl_2}^{\!*}}^{\,\Z} \Big( \!
\Big(\! {{\Lambda_2 \, ; \, 0} \atop n} \Big) \Big\vert_{q=\varepsilon}
\Big) = \, \cases
     \hskip-3pt  \Big(\! {{\Lambda_2 \, ; \; 0} \atop {n / \ell}} \!\Big)
{\Big|}_{q=1}  \hskip-3pt  = \Big(\! {\g_2 \atop {n / \ell}}
\!\Big)  &  \text{if} \quad \ell \,\big\vert\, n  \\
     \hskip1pt  \quad\;  0  &  \text{if} \quad \ell \!\not\big\vert\, n
   \endcases  $$
 \vskip-2pt
\noindent
%
%
 while  $ \, {\widehat{\xi}\,{\Big|}_{\! q=1}}\! \! \Big(\! \Big(\!
{{d \, ; \, 0} \atop {n / \ell}} \Big) \Big\vert_{q=1} \Big)\! =\!
\Big(\! {{\Lambda_2 \, ; \, 0} \atop {n / \ell}} \Big) \Big\vert_{q=1}
\hskip-5pt = \, \Big(\! {\g_2 \atop {n / \ell}} \!\Big) \, $  due to
Corollary 3.6{\it (a)},  so  $ \, {\Cal{F}r}_{M_2}^{\,\Z}\! \Big(\!
\Big(\! {{d \, ; \, 0} \atop n} \Big)\Big\vert_{q=\varepsilon} \Big)
= \Big(\! {{d \, ; \, 0} \atop {n / \ell}} \Big) \Big\vert_{q=1}
\hskip-3pt = \Big(\! {\g_2 \atop {n / \ell}} \!\Big) \; $  whenever
$ \, \ell \big| n \, $,  \, and otherwise  $ \; {\Cal{F}r}_{M_2}^{\,\Z}
\Big( \! \Big(\! {{d \, ; \, 0} \atop n} \Big) \Big\vert_{q=\varepsilon}
\Big) = \, 0 \, $.
                                             \par
   All this accounts for claim  {\it (a)}.  Claims  {\it (b)\/}  and
{\it (c)\/}  can be proved with the same arguments, or deduced from
{\it (a)\/}  in force of Proposition 3.7 and of Corollary 3.3.   \qed
\enddemo

\vskip1,1truecm

\centerline {\bf \S\; 4 \ Miscellanea results on  $ q $--numbers
and  $ q $--functions.}

\vskip13pt

   In several steps along the present work we need special, technical
results about  $ q $--numbers  and their combinatorics.  We collect them
in the present section, referring to [GR1] for proofs.

\vskip11pt

 \proclaim{Lemma 4.1}  For all  $ \, k \in \N \, $,  \, let  $ \,
\Pi_k := {(q-1)}^k {(k)}_q! \, $  and  $ \, (x \, ; k) := \Pi_k
\cdot \Big(\! {{x \, ; \, 0} \atop k} \!\Big) \, $.  Then
  $$  (x \, ; n) \; = \; {\textstyle \sum\limits_{k=0}^n} \; {(-1)}^{n-k}
\, q^{-{k \choose 2}} \, \textstyle{\Big( {n \atop k} \Big)}_{\!q^{-1}}
\, x^k  \quad ,  \qquad
   x^n \; = \; {\textstyle \sum\limits_{k=0}^n} \; q^{k \choose 2} \,
\textstyle{\Big( {n \atop k} \Big)}_{\!q} \, (x \, ; \,k)  \quad .
\qed  $$
\endproclaim

\vskip11pt

 \proclaim{Lemma 4.2}  Let  $ A $  be any  $ \Qq $--algebra,  and let
$ \, x, y, z, w \in A \, $  be such that  $ \; x \, w = q^2 \, w \,
x \, $,  $ \; x \, y = q \, y \, x \, $,  $ \; x \, z = q \, z \, x \; $
and  $ \; y \, z = z \, y \, $.  Then for all  $ \, n \in \N \, $  and
$ \, t \in \Z \, $  we have
  $$  \leqalignno{
   \textstyle{ \Big( {{x \, + \, {( q - q^{-1} \,)}^2 \, w \; ; \; t} \atop n}
\Big)} \;\; =  &  \;\; {\textstyle \sum\limits_{r=0}^n} \; q^{r(t-n)}
\big( q - q^{-1} \big)^r \cdot \, w^{(r)} \, \textstyle{\Big\{
{{x \, ; \; t} \atop {n \, , \, r}} \Big\}} \;\; =   &  \hbox{\it (a-1)}  \cr
   =  &  \;\; {\textstyle \sum\limits_{r=0}^n} \; q^{r(t-n)} \big( q - q^{-1}
\big)^r \cdot \, \textstyle{\Big\{ {{x \, ; \; t - 2 \, r} \atop {n \, , \, r}}
\Big\}} \; w^{(r)}   &   \hbox{\it (a-2)}  \cr
   \textstyle{ \Big( {{x \, + \, {( q - q^{-1} \,)}^2 \, y \, z \; ; \; t} \atop n}
\Big)} \;\; =  &  \;\; {\textstyle \sum\limits_{r=0}^n} \; q^{r(t-n)}
\big( q - q^{-1} \big)^r \, {[r]}_q! \cdot y^{(r)} \,
{\textstyle \Big\{ {{x \, ; \; t - r } \atop {n \, , \, r}} \Big\}}
\; z^{(r)}   \qquad .  \qed   &  \hbox{\it (b)}  \cr }  $$
\endproclaim

\vskip7pt

 \proclaim{Lemma 4.3}  Let  $ \varOmega $  be any  $ \Zqqm $--algebra,
$ \, \varepsilon \, $  be a (formal) primitive  $ \ell $-th  root of 1,
with  $ \, \ell \in \N_+ \, $,  \, and  $ \, \varOmega_\varepsilon :=
\varOmega \Big/ (q - \varepsilon) \, \varOmega \, $  the specialization
of  $ \varOmega $  at  $ \, q = \varepsilon \, $.  Then for each  $ \,
x, y \in \varOmega \, $
  $$  \varOmega \ni {\textstyle \Big(\! {{x \, ; \, 0} \atop \ell} \!\Big)}
\;\; \Longrightarrow \;\, \big(x\big|_{q=\varepsilon}\big)^\ell = \, 1
\text{\ \ in \ } \! \varOmega_\varepsilon \, ,  \qquad \!  \varOmega \ni
y^{(\ell)} \;\; \Longrightarrow \;\, \big( y\big|_{q=\varepsilon} \big)^\ell
= \, 0 \text{\ \ in \ } \! \varOmega_\varepsilon \; .   \qed  $$
\endproclaim

\vskip3pt

\proclaim{Lemma 4.4}  For all  $ \, n \in \N \, $,  \, the identity
$ \;\; a^n \, d^n = \sum\limits_{k=0}^n \, q^{k^2}
\Big(\! {n \atop k} \!\Big)_{\!q^2} \, b^k \, c^k \;\; $
holds in  $ \, \fqsl \, $.   $ \square $
\endproclaim

\vskip11pt

\proclaim{Lemma 4.5}  The following identities holds (notation of \, \S\; 1.2):
  $$  q^j \, {\textstyle \Big\{\! {{x \, ; \; s} \atop {m \, , \, j}} \!\Big\}}
\; \big( q^{s-m+2j} \, x - 1 \big) \, + \, \big( q^{2j} - 1 \big) \,
{\textstyle \Big\{\! {{x \; ; \; s} \atop {m \, , \, j-1}} \!\Big\}}
\; = \; \big( q^{m+1} - 1 \big) \, {\textstyle \Big\{ {{x \; ; \; s}
\atop {m+1 \, , \, j}} \Big\}}   \leqno \text{\it (a)}  $$
 \vskip-11pt
  $$  q^{s \choose 2} \, {\textstyle \Big(\! {n \atop s} \!\Big)}_{\!q}
\;  =  \; {\textstyle \sum\limits_{j=1}^s} \; {(-1)}^{j-1} \, q^{{s-j}
\choose 2} \, {\textstyle \Big(\! {n \atop j} \!\Big)}_{\!q} \,
{\textstyle \Big(\! {{n-j\,} \atop {s-j}} \!\Big)}_{\!q}
\leqno \text{\it (b)}  $$
 \vskip-11pt
  $$  {\textstyle \sum\limits_{j=0}^k} \, {(-1)}^j \, q^{j \choose 2}
\, {\textstyle \Big(\! {{n-j\,} \atop {k-j}} \!\Big)}_{\!q^2} \,
{\textstyle \Big(\! {n \atop j} \!\Big)}_{\!q} \;  =  \;\,
q^{k^2} \, {(q-1)}^k \, {\textstyle \Big(\! {n \atop {2\,k}} \!\Big)}_{\!q}
{\big(2k-1\big)}_q!!   \hskip25pt \forall \;\; k \leq n
\leqno \text{\it (c)}  $$
\endproclaim

\demo{Proof}  To give an idea, we sketch the proof of  {\it (c)\/}   --- the
complete proof is in [GR1].  First we fix some more notation: for all  $ \, s,
k, n \in \N \, $  set  $ \; {\langle s \rangle}_q := {{\, q^s + 1 \,} \over
{\, q + 1 \,}} \; $,  $ \;\; {\langle s \rangle}_q! := \prod_{r=1}^s
{\langle r \rangle}_q \; $,  $ \;\; \Big\langle\! {n \atop k}
\!\Big\rangle_{\!q} := {{\, {\langle n \rangle}_q! \,} \over
{\, {\langle k \rangle}_q! \, {\langle n\!-\!k \rangle}_q! \,}}
\; ,  \;\; {(x \, ; n)}_q := \prod_{s=1}^n \big( x \, q^{1-s}
\! - 1 \big) \, $,  $ \;\; {\langle x \, ; n \rangle}_q :=
\prod_{s=1}^n \big( x \, q^{1-s} + 1 \big) \; $.  The following
identities then are clear from definitions:
  $$  \displaylines{
 \textstyle
   {\Big(\! {n \atop k} \!\Big)}_{\!q^2} = \,
{\Big(\! {n \atop k} \!\Big)}_{\!q} \,
{\Big\langle {n \atop k} \Big\rangle}_{\!q} \;\; ,
\quad  {\Big(\! {n \atop j} \!\Big)}_{\!q} \,
{\Big(\! {{n-j\,} \atop {k-j\,}} \!\Big)}_{\!q} = \,
{\Big(\! {n \atop k} \!\Big)}_{\!q} \,
{\Big(\! {k \atop j} \!\Big)}_{\!q} \;\; ,
\quad  {\Big(\! {n \atop k} \!\Big)}_{\!q} =
{{\,{( q^n ; \, k )}_q\,} \over {\,{( q^k ; \, k )}_q\,}}  \cr
 \textstyle
   \prod_n^- := \prod_{i=1}^n
\big( q^j - 1 \big) = {\big( q^n ; n \big)}_q \;\; ,  \quad
\prod_n^+ := \prod_{i=1}^n
\big( q^j + 1 \big) = {\big\langle q^n ; n \big\rangle}_q \;\; ,  \quad
{\Big\langle\! {n \atop k} \!\Big\rangle}_{\!q} = {{\,{\langle q^n \, ;
\, k \,\rangle}_q\,} \over {\,{\langle q^k ; \, k \,\rangle}_q\,}}
\cr }  $$
   \indent   Using them, we transform a bit both sides of  {\it (c)}.  Namely,
we get
  $$  \displaylines{
   \qquad  {\textstyle \sum_{j=0}^k} \, {(-1)}^j \, q^{j \choose 2}
\, {\textstyle {\Big( {{n-j\,} \atop {k-j\,}} \!\Big)}}_{\!q^2} \,
{\textstyle {\Big(\! {n \atop j} \!\Big)}}_{\!q}  =  \;
{\textstyle \sum_{j=0}^k} \, {(-1)}^j \, q^{j \choose 2} \,
{\textstyle {\Big\langle\! {{n-j} \atop {k-j}} \Big\rangle}_{\!q}}
\, {\textstyle {\Big(\! {{n-j\,} \atop {k-j\,}} \!\Big)}_{\!q}} \,
{\textstyle {\Big(\! {n \atop j} \!\Big)}_{\!q}} \,  =   \hfill  \cr
   \hfill   =  \, {\textstyle \sum_{j=0}^k} \, {(-1)}^j \,
q^{j \choose 2} \, {\textstyle {\Big(\! {n \atop k} \!\Big)}_{\!q}}
\, {\textstyle {\Big( {k \atop j} \Big)}_{\!q}} \,
{\textstyle {\Big\langle\! {{n-j} \atop {k-j}} \Big\rangle}_{\!q}}
\,  =  \; {\textstyle {\Big(\! {n \atop k} \!\Big)}_{\!q}} \;
{\textstyle \sum_{j=0}^k} \, {(-1)}^j \, q^{j \choose 2} \,
{\textstyle {\Big(\! {k \atop j} \!\Big)}_{\!q}} \,
{\textstyle {\Big\langle\! {{n-j} \atop {k-j}} \Big\rangle}_{\!q}} }  $$
for the l.h.s., and similarly for the r.h.s.~we find
                                             \par
   \centerline{ $ q^{k^2} \, {(q-1)}^k \, x \, {\Big(\! {n \atop {2\,k}}
\!\Big)}_{\!q} \, {\big( 2k \! - \!1\big)}_q!! \;  =  \; q^{k^2} \, {(q-1)}^k
\, {{\, {(n)}_q \, {(n \! - \! 1)}_q \cdots {(n-2k+1)}_q \,} \over
{\,{(2\,k)}_q!!\,}} \;  =  \; q^{k^2} \, {\Big( {n \atop k} \Big)}_{\!q}
\, {\Big(\! {{n-k} \atop k} \!\Big)}_{\!q} \cdot {{\, \prod_k^- \,} \over
{\, \prod_k^+ \,}} $ }
\noindent
 In light of this, we must prove that  $ \; \sum_{j=0}^k \, {(-1)}^j \,
q^{j \choose 2} \, \Big(\! {k \atop j} \!\Big)_{\!q} \, {\Big\langle\!
{{n-j} \atop {k-j}} \Big\rangle}_{\!q} \,  =  \, q^{k^2} \, \Big(\!
{{n-k} \atop k} \!\Big)_{\!q} \cdot {{\, \prod_k^- \,} \over
{\, \prod_k^+ \,}} \; $,
\, or
  $$  {\textstyle \sum\limits_{j=0}^k} \, {(-1)}^{k-j} \, q^{{k-j} \choose 2}
\, {\textstyle {\Big(\! {k \atop j} \!\Big)}_{\!q}} \, {\big\langle q^{n-k+j}
\, ; j \,\big\rangle}_q \, {\big\langle q^k \, ; k-j \, \big\rangle}_q \;  =
\;\, q^{k^2} \, {\big( q^{n-k} \, ; k \,\big)}_q  \quad .   \eqno (4.1)  $$
To this end, we shall prove a more general identity.  Let
  $$  R^n_{k,t} \;  :=  \; {\textstyle \sum_{j=0}^k} \, {(-1)}^{k-j} \,
q^{{k-j} \choose 2} \, {\textstyle {\Big(\! {k \atop j} \!\Big)}_{\!q}}
\, {\big\langle q^{n-k+j} \, ; t+j \,\big\rangle}_q \, {\big\langle
q^{k+t} \, ; k-j \,\big\rangle}_q  \quad ;   \eqno (4.2)  $$
we shall prove that
  $$  R^n_{k,t} \;  =  \; q^{k(k+t)} \, {\big\langle q^{n-k} \, ;
\, t \big\rangle}_q \; {\big( q^{n-k-t} \, ; \, k \big)}_q  \quad .
\eqno (4.3)  $$
 \eject
\noindent
 Note that, by definition, (4.3) for  $ \, t = 0 \, $  yields exactly (4.1),
so the latter is just a special case of the former.  The proof of (4.3)
follows easily by induction on  $ k $  from the
following identity:
  $$  \qquad \hskip25pt   R_{(k+1),t}^n \;  =  \; - \, q^k \,
\big( q^{k+t+1} + 1 \big) \, R_{k,t}^{n-1} \, + \, R_{k,t+1}^n   \qquad
(\, \forall \; k, t, n \,) \; .   \eqno (4.4)  $$
   \indent   Let us prove now the above identity.  We have
  $$  \displaylines{
 \textstyle
   R_{(k+1),t}^n \;  =  \; \sum_{j=0}^{k+1} \, {(-1)}^{k+1-j} \,
q^{{{k+1-j} \choose 2}} \, \Big(\! {{k+1} \atop j} \!\Big)_{\!q}
\, {\big\langle q^{n-k-1+j} \, ; \, t+j \,\big\rangle}_q \;
{\big\langle q^{k+1+t} \, ; \, k+1-j \,\big\rangle}_q \;  =   \hfill  \cr
 \textstyle
   \hfill   =  \; \sum\limits_{j=0}^{k+1} \, {(-1)}^{k+1-j} \,
q^{{{k+1-j} \choose 2}} \, \Big(\! \Big( {k \atop {j-1}} \Big)_{\!q}
+ \, q^j \, \Big( {k \atop j} \Big)_{\!q} \, \Big) \,
{\big\langle q^{n-k-1+j} \, ; \, t+\!j \,\big\rangle}_q
\, {\big\langle q^{k+1+t} \, ; \, k+1-\!j \,\big\rangle}_q
\,  =  \quad  \cr
 \textstyle
   \quad   =  \; - \, q^k \, \sum_{j=0}^k \, {(-1)}^{k-j} \, q^{{{k-j} \choose 2}}
\, \Big( {k \atop j} \Big)_{\!q} \, {\big\langle q^{n-k-1+j} \, ; \, t+j
\,\big\rangle}_q \, {\big\langle q^{k+1+t} \, ; \, k+1-j \,\big\rangle}_q
\, +   \hfill  \cr
 \textstyle
   \hfill   + \, \sum_{j=1}^{k+1} \, {(-1)}^{k+1-j} \, q^{{{k+1-j} \choose 2}}
\, \Big( {k \atop {j-1}} \Big)_{\!q} \, {\big\langle q^{n-k-1+j} \, ; \, t+j
\,\big\rangle}_q \, {\big\langle q^{k+1+t} \, ; \, k+1-j \,\big\rangle}_q
\;  =  \quad  \cr
 \textstyle
   \quad   =  \; - \, q^k \, \big( q^{1+k+t} + 1 \big) \, \sum_{j=0}^k \,
{(-1)}^{k-j} \, q^{{{k-j} \choose 2}} \, \Big( {k \atop j} \Big)_{\!q}
\, {\big\langle q^{(n-1)-k+j} \, ; \, t+j \,\big\rangle}_q \,
{\big\langle q^{k+t} \, ; \, k-j \,\big\rangle}_q \, +   \hfill  \cr
 \textstyle
   \hfill   + \, \sum_{j=0}^k \, {(-1)}^{k-j} \, q^{{{k-j} \choose 2}} \,
\Big( {k \atop j} \Big)_{\!q} \, {\big\langle q^{n-k+j} \, ; \, (t+1)+j
\,\big\rangle}_q \, {\big\langle q^{k+(t+1)} \, ; \, k-j \,\big\rangle}_q
\;  =  \quad  \cr
   \hfill   =  \; - \, q^k \, \big( q^{k+t+1} + 1 \big) \, R_{k,t}^{n-1}
\, + \, R_{k,t+1}^n  \;\; ,  \quad \text{q.e.d.} }  $$
   \indent   Now the induction on  $ k $  to prove (4.3) goes as follows.
For all  $ \, n \in \N_+ \, $  and all  $ \, t \in \N \, $,  the right
hand sides of (4.2) and (4.3) are equal, because in that case (4.4)
gives
  $$  \displaylines{
   - \, \big( q^{n-1} + 1 \big) \cdots \big( q^{nt+1} + 1 \big) \big(
q^{t+1} + 1 \big) + \big( q^{n} + 1 \big) \cdots \big( q^{n-t} + 1 \big)
\, =   \hfill  \cr
  \hfill   = \, \big( q^{n-1} + 1 \big) \cdots \big( q^{n-t} + 1 \big) \,
\big( q^{n} + 1 - q^{t+1} - 1 \big) \, = \, q^{t+1} \, \big( q^{n-1} +
1 \big) \cdots \big( q^{n-t} + 1 \big) \, \big( q^{n-t-1} - 1 \big)
\; .  \cr }  $$
This sets the basis of induction.  For the inductive step, assume the
r.h.s.'s of (4.2) and (4.3) be equal for all  $ \, t \in \N \, $  and
$ \, k', n \in \N_+ $  with  $ \, k' \leq k \, $,  \, for some  $ \,
k \in \N_+ \, $,  $ \, k \leq n \, $.  Then
  $$  \displaylines{
   R_{(k+1),t}^n \;  =  \; - \, q^k \, \big( q^{k+t+1} + 1 \big) \,
R_{k,t}^{n-1} + R_{k,t+1}^n \;  =   \hfill  \cr
   =  \; - \, q^k \, \big( q^{k+t+1} + 1 \big) \, q^{k(k+t)} \,
{\big\langle q^{n-1-k} \, ; t \big\rangle}_q \, {\big( q^{n-1-k-t} \, ; k \big)}_q
\, +  \cr
   \hfill   + \; q^{k(k+(t+1))} \, {\big\langle q^{n-k} ; \, t+1 \big\rangle}_q
\, {\big( q^{n-k-t-1} \, ; k \big)}_q \;  =  \cr
   \, =  - \, q^{k(k+t+1)} \, \Big( \big( q^{k+t+1} + 1 \big) \,
{\big\langle q^{n-1-k} \, ; t \big\rangle}_q \,
{\big( q^{n-1-k-t} \, ; k \big)}_q \, + \,
{\big\langle q^{n-k} \, ; t+1 \big\rangle}_q \,
{\big( q^{n-k-t-1} \, ; k \big)}_q \Big)  =  \cr
   \qquad  =  \; q^{k(k+(t+1))} \, {\big\langle q^{n-1-k} \, ; t \big\rangle}_q
\, {\big( q^{n-1-k-t} \, ; k \big)}_q \, \big( q^{n-k} + 1 - 1 - q^{k+t+1} \big)
\;  =   \hfill  \cr
   =  \; q^{k^2 + k \, (t+1) + k + t + 1} \,
{\big\langle q^{n-1-k} \, ; t \big\rangle}_q \,
{\big( q^{n-1-k-t} \, ; k \big)}_q \, \big( q^{n - 2(k+1) - t + 1} - 1 \big)
\;  =  \cr
   \hfill   =  \; q^{(k+1) \, ((k+1)+t)} \,
{\big\langle q^{n-(k+1)} \, ; t \big\rangle}_q
\, {\big( q^{n-(k+1)-t} \, ; k+1 \big)}_q }  $$
which gives exactly (4.3) with  $ \, k+1 \, $  instead of  $ k \; $.   \qed
\enddemo

\vskip1,4truecm

\vskip37pt

\enddocument

\Refs
  \widestnumber\key {GR2}

\vskip3pt

\ref
 \key  CP   \by  V. Chari, A. Pressley
 \book  A guide to Quantum Groups
 \publ  Cambridge Univ. Press
 \publaddr  Cambridge   \yr  1994
\endref

\vskip1pt

\ref
 \key  DL   \by  C. De Concini, V. Lyubashenko
 \paper  Quantum function algebra at roots of 1
 \jour  Adv. Math.   \vol  108   \yr  1994
 \pages  205--262
\endref

\vskip1pt

\ref
 \key  DP   \by  C. De Concini, C. Procesi
 \paper  Quantum groups
 \inbook  D-modules, Representation Theory, and Quantum Groups
 \eds  L. Boutet de Monvel, C. De Concini, C. Procesi, P. Schapira, M. Vergne
 \yr  1993   \pages  798--820
 \publ  Lecture Notes in Mathematics  {\bf 1565},  Springer \& Verlag
 \publaddr  Berlin
%
\endref

\vskip1pt

\ref
 \key  Dr   \by  V. G. Drinfel'd
 \paper  Quantum groups
 \inbook  Proc. Intern. Congress of Math. (Berkeley, 1986)  \yr  1987
 \pages  798--820
\endref

\vskip1pt

\ref
 \key  Ga1   \by  F. Gavarini
 \paper  Quantization of Poisson groups
 \jour  Pacific Journal of Mathematics
 \vol  186   \yr  1998   \pages  217--266
\endref

\vskip1pt

\ref
 \key  Ga2   \bysame   
 \paper  Quantum function algebras as quantum enveloping algebras
 \jour  Communications in Algebra
 \vol  26   \yr  1998   \pages  1795--1818
\endref

\vskip1pt

\ref
 \key  Ga3   \bysame   
 \paper  The global quantum duality principle: theory,
examples, and applications
 \jour  preprint\break   {\tt math.QA/0303019}   \yr  2003
\endref

\vskip1pt

\ref
 \key  Ga4   \bysame   
 \paper  The global quantum duality principle
 \jour   J. Reine Angew. Math.   \toappear
\endref

\vskip1pt

\ref
 \key  GR1   \by  F. Gavarini, Z. Raki\'c
 \paper  $ F_q[{M}_2] $,  $ F_q[{GL}_2] $  and
$ F_q[{SL}_2] $  as quantized hyperalgebras
 \jour  preprint\break  {\tt math.QA/0411440}
 \yr  2004
\endref

\vskip1pt

\ref
 \key  GR2   \by  F. Gavarini, Z. Raki\'c
 \paper  $ F_q[{M}_n] $,  $ F_q[{GL}_n] $  and
$ F_q[{SL}_n] $  as quantized hyperalgebras
 \jour  J. Algebra   
 \vol  315   
 \yr  2007   
 \pages  761--800   
\endref

\vskip1pt

\ref
 \key  Hu   \by  J. E. Humphreys
 \book  Introduction to Lie Algebras and Representation Theory
 \publ  Graduate Texts in Mathematics {\bf 9},  Springer \& Verlag
 \publaddr  Berlin--Heidelberg--New York   \yr  1972
\endref

\vskip1pt

\ref
 \key  No   \by  M. Noumi
 \paper  Macdonald's Symmetric Polynomials as Zonal Spherical
Functions on Some Quantum Homogeneous Spaces
 \jour  Adv. Math.   \vol  123   \yr  1996   \pages  16--77
\endref

\endRefs
\bye
\end